\documentclass[12pt]{amsart}
\date{22 November 2006}
\usepackage{latexsym,amsmath,amsfonts,amscd,amssymb}
\usepackage{lscape}
\usepackage{mathdots}
\usepackage{array}
\DeclareFontFamily{OT1}{rsfs}{} \DeclareFontShape{OT1}{rsfs}{n}{it}{<->rsfs10}{}
 \DeclareMathAlphabet{\curly}{OT1}{rsfs}{n}{it}

\theoremstyle{plain}  
\newtheorem{theorem}{Theorem}[section]

\newtheorem*{theorem*}{Theorem}

\newtheorem{corollary}[theorem]{Corollary}

\newtheorem{proposition}[theorem]{Proposition}

\theoremstyle{remark}

\newtheorem{remark}[theorem]{Remark}

\newtheorem*{claim*}{Claim}
\numberwithin{equation}{section}
\newcommand{\suchthat}{\;:\;}
\newcommand{\abs}[1]{\lvert#1\rvert}
\renewcommand{\leq}{\leqslant}

\renewcommand{\geq}{\geqslant}

\newcommand{\R}{\mathbb{R}}

\newcommand{\Z}{\mathbb{Z}}
\newcommand{\C}{\mathbb{C}}
\newcommand{\E}{\mathbb{E}}
\newcommand{\GGG}{\curly{G}}
\newcommand{\HHH}{\curly{H}}
\newcommand{\HH}{\mathbb{H}}

\newcommand{\cD}{\mathcal{D}}
\newcommand{\cM}{\mathcal{M}}
\newcommand{\cO}{\mathcal{O}}
\newcommand{\cR}{\mathcal{R}}

\newcommand{\Sh}{\check{S}}
\newcommand{\bE}{{\bf{E}}}
\newcommand{\dbar}{\bar{\partial}}
\newcommand{\lie}{\mathfrak}
\newcommand{\ra}{\rightarrow}
\newcommand{\lra}{\longrightarrow}

\newcommand{\PSL}{\mathrm{PSL}}

\newcommand{\SU}{\mathrm{SU}}
\newcommand{\U}{\mathrm{U}}
\newcommand{\OO}{\mathrm{O}}
\newcommand{\GL}{\mathrm{GL}}
\newcommand{\SL}{\mathrm{SL}}
\newcommand{\SSS}{\mathrm{S}}
\newcommand{\SO}{\mathrm{SO}}
\newcommand{\Sp}{\mathrm{Sp}}
\newcommand{\Spin}{\mathrm{Spin}}

\DeclareMathOperator{\ad}{ad}

\DeclareMathOperator{\rk}{rk}
\DeclareMathOperator{\rank}{rank}
\DeclareMathOperator{\im}{im}

\DeclareMathOperator{\Hom}{Hom}
\DeclareMathOperator{\End}{End}

\DeclareMathOperator{\Herm}{Herm}
\DeclareMathOperator{\Mat}{Mat}
\DeclareMathOperator{\Sym}{Sym}
\DeclareMathOperator{\Skew}{Skew}
\newcommand{\norm}[1]{\lVert#1\rVert}

\renewcommand{\phi}{\varphi}

\newcommand{\liem}{\mathfrak{m}}
\newcommand{\liemp}{\mathfrak{m}_+}
\newcommand{\liemm}{\mathfrak{m}_-}
\newcommand{\liemc}{\mathfrak{m}^{\mathbb{C}}}
\newcommand{\lieh}{\mathfrak{h}}
\newcommand{\liehc}{\mathfrak{h}^{\mathbb{C}}}
\newcommand{\lieg}{\mathfrak{g}}
\newcommand{\liegc}{\mathfrak{g}^{\mathbb{C}}}
\newcommand{\vol}{\mathrm{vol}}

\begin{document}

\title[Maximal surface group representations]
{Maximal surface group representations in isometry groups of
classical Hermitian symmetric spaces}

\author[S. B. Bradlow]{Steven B. Bradlow}
\address{Department of Mathematics \\
University of Illinois \\
Urbana \\
IL 61801 \\
USA }
\email{bradlow@math.uiuc.edu}

\author[O. Garc{\'\i}a-Prada]{Oscar Garc{\'\i}a-Prada}
\address{Departamento  de Matem\'aticas \\
  CSIC \\ Serrano 121 \\ 28006 Madrid \\ Spain}
\email{oscar.garcia-prada@uam.es}

\author[P. B. Gothen]{Peter B. Gothen}
\address{Departamento de Matem\'atica Pura \\
  Faculdade de Ci\^encias, Universidade do Porto \\
  Rua do Campo Alegre 687 \\ 4169-007 Porto \\ Portugal }
\email{pbgothen@fc.up.pt}

\thanks{
  Members of VBAC (Vector Bundles on Algebraic Curves).
  Second and Third authors partially supported by Ministerio de
  Educaci\'{o}n y Ciencia and Conselho de Reitores das
  Universidades Portuguesas through Acci\'{o}n Integrada Hispano-Lusa
  HP2002-0017 (Spain) / E--30/03 (Portugal).
  First and Second authors partially supported by Ministerio de Educaci{\'o}n
y  Ciencia (Spain) through Project MTM2004-07090-C03-01.
  Third author partially supported by the Centro de Matem\'atica da
  Universidade do Porto and the project POCTI/MAT/58549/2004, financed
  by FCT (Portugal) through the programmes POCTI and POSI of the QCA
  III (2000--2006) with European Community (FEDER) and national
  funds.
The second author visited the IHES with the partial support of the European
Commission through its 6th Framework Programme ``Structuring the
European Research Area'' and the Contract No.\ RITA-CT-2004-505493 for
the provision of Transnational Access implemented as Specific Support
Action.
}


\begin{abstract}Higgs bundles and non-abelian Hodge theory provide
holomorphic methods with which to study the moduli spaces of surface group
representations in a reductive Lie group $G$. In this paper we survey
the case in which $G$ is the isometry group of a classical Hermitian 
symmetric space of non-compact type. Using Morse theory on the moduli 
spaces of Higgs bundles, we compute the number of connected
components of the moduli space of representations with maximal Toledo
invariant. 
\end{abstract}

\maketitle

\section{Introduction}

Given a closed oriented surface, $X$, and a connected semisimple  Lie group
$G$, the moduli space of representations of $\pi_1(X)$ in $G$ is defined as
the set
$$
\mathcal{R}(G) = \Hom^{+}(\pi_1(X),G) / G
$$
of reductive homomorphisms from $\pi_1(X)$ to $G$ modulo conjugation. The
reductiveness condition ensures that this orbit space is Hausdorff, and
in fact $\mathcal{R}(G)$ is a real analytic variety. The geometry and
topology of these moduli spaces, though clearly reflective of properties of
both $X$ and $G$, is still far from fully understood, especially in the
case where $G$ is non-compact. In this paper we consider the non-compact
groups for which the homogeneous space $G/H$, where $H\subset G$ is a
maximal compact subgroup, is a Hermitian symmetric space.  By the Cartan
classification of irreducible symmetric spaces, it thus suffices for us to
consider the groups $G=\SU(p,q)$, $G=\Sp(2n,\R)$, $G=\SO^*(2n)$ and
$G=\SO_0(2,n)$. We concentrate mainly on the most primitive topological
property of $\mathcal{R}(G)$, namely the number of connected components.

The first division of $\mathcal{R}(G)$ into disjoint closed subspaces
comes from the correspondence between representations of $\pi_1(X)$
and flat principal bundles over $X$. Every representation $\rho:
\pi_1(X)\ra G$ carries a topological invariant which is the
characteristic class $d\in\pi_1(G)$ of the flat $G$-bundle
corresponding to $\rho$. This class measures the obstruction to lift
$\rho$ to a representation of $\pi_1(X)$ in the universal cover of
$G$. The subvarieties $\mathcal{R}_d(G)\subset \mathcal{R}(G)$,
consisting of representations with a fixed value of the invariant,
form disjoint closed subspaces but not necessarily connected
components. The problem is to count and understand the distinct
components of each of the $\mathcal{R}_d(G)$.

For compact groups $G$, it is well-known that for every $d\in\pi_1(G)$, the
moduli space $\mathcal{R}_d(G)$ is non-empty and connected. One way to see
this is to choose a complex structure on $X$ and to use the theory of
holomorphic bundles on the resulting Riemann surface. When $G=\SU(n)$, the
topological invariant is trivial, since $\SU(n)$ is simply connected, and
by a theorem of Narasimhan and Seshadri \cite{narasimhan-seshadri:1965},
$\mathcal{R}(\SU(n))$ can be identified with the moduli space of polystable
vector bundles of rank $n$ and trivial determinant, which is connected. A
similar result was proved by Ramanathan \cite{ramanathan:1975} for every
connected compact semisimple Lie group $G$. He identified
$\mathcal{R}_d(G)$ with the moduli space of polystable holomorphic
principal $G^\C$-bundles over $X$ with topological class $d\in \pi_1(G)$,
where $G^\C$ is the complexification of $G$, and showed that this moduli
space is connected.

Holomorphic methods can also be used when $G$ is a complex semisimple
Lie group. In place of a holomorphic principal bundle, the holomorphic object
corresponding to a representation is now a $G$-Higgs bundle, i.e.\ a
pair consisting of a holomorphic $G$-bundle and a holomorphic
section of the adjoint bundle twisted with the canonical bundle of
$X$. A combination of theorems by Hitchin \cite{hitchin:1987a} and
Donaldson \cite{donaldson:1987} for $G=\SL(2,\C)$ and Simpson
\cite{simpson:1992} and Corlette \cite{corlette:1988} for general
$G$ identify $\mathcal{R}_d(G)$ with $\mathcal{M}_d(G)$, the moduli
space of polystable $G$-Higgs bundles with fixed topological class.
Morse-theoretic methods introduced by Hitchin \cite{hitchin:1987a}
prove the connectedness of $\mathcal{M}_d(G)$ by relating it to that
of the moduli space of polystable holomorphic principal
$G$-bundles\footnotemark. \footnotetext{Other methods to prove that
$\mathcal{M}_d(G)$ is
  connected had been used by Goldman \cite{goldman:1988} for
  $G=\SL(2,\C)$ and $G=\PSL(2,\C)$ and by J. Li \cite{jun-li:1993} for
  an arbitrary semisimple complex Lie group $G$.}

The situation is very different if $G$ is a non-compact real form of a
semisimple complex Lie group. The simplest case is $G=\SL(2,\R)$. In this
case $\pi_1(G)=\Z$ and the topological invariant $d\in\Z$ of a
representation is the Euler class of the corresponding flat
$\SL(2,\R)$-bundle. The Milnor--Wood inequality says that
$\mathcal{R}_d(G)$ is empty unless $|d|\leq g-1$, where $g$ is the genus of
$X$. In \cite{goldman:1980}, Goldman showed that
$\mathcal{R}_{\pm(g-1)}(G)$ has $2^{2g}$ components consisting of discrete
faithful representations, each of which can be identified with the
Teichm\"uller space of $X$. Later in \cite{goldman:1988} he showed that
$\mathcal{R}_d(G)$ is connected for $|d|< g-1$. This was also proved by
Hitchin \cite{hitchin:1987a} using Higgs bundle methods.

The results for $G=\SL(2,\R)$ can be generalized in two ways. The
first goes back to \cite{hitchin:1992} where Hitchin extended the
results from $\SL(2,\R)$ to $G=\SL(n,\R)$.  Using Higgs bundles, he
counted the number of connected components and, moreover, for any
split real form identified a component homeomorphic to $\R^{\dim G
  (2g-2)}$ and which naturally contains a copy of Teichm\"uller space.
This component, known as the Teichm\"uller or Hitchin component, has
special geometric significance. In the case of $\SL(3,\R)$, Choi and
Goldman \cite{goldman-choi:1993,goldman-choi:1997} showed that the
representations in the Hitchin component are discrete and faithful and
correspond to convex projective structures on the surface. More
recently, Labourie \cite{labourie:2004} has shown for $G=\SL(n,\R)$
that these representations are discrete and faithful and are related
to certain Anosov geometric structures, and analogous results have
been obtained for $G = \Sp(2n,\R)$ by Burger, Iozzi, Labourie and
Wienhard \cite{burger-iozzi-labourie-wienhard:2005}.

The second generalization from $G=\SL(2,\R)$ exploits the fact that
the homogeneous space $G/H$, where $H\subset G$ is a maximal compact
subgroup, is Hermitian symmetric. Indeed, $\SL(2,\R)/\SO(2)$ is the
hyperbolic plane. This generalization is the main focus of this paper.
In particular, we review how the theory of Higgs bundles is used in
the general case of a non-compact real form $G$ such that the
symmetric space $G/H$ is Hermitian.

We start by introducing the appropriate notion of a $G$-Higgs bundle when
$G$ is any connected reductive real Lie group. The information required to
define such a Higgs bundle includes a choice of maximal compact subgroup
$H\subset G$ and a Cartan decomposition $\lieg=\lieh +\liem$ of the Lie
algebras.  The correct notion of a Higgs bundle then turns out to be a pair
$(E,\varphi)$, where $E$ is an $H^\C$-bundle (where $H^\C$ is the
complexification of $H$) and the Higgs field $\varphi$ takes values in the
complexification $\liem^\C$ of $\liem$.

Topological classes of $E$ are characterized by elements $d\in
\pi_1(H)\cong \pi_1(G)$. When $G$ is semisimple and $G/H$ is an
irreducible Hermitian symmetric space, the torsion-free  part of
$\pi_1(H)$ is isomorphic to $\Z$ (in fact, for all the classical
groups that we will study $\pi_1(H)\cong \Z$, except in the cases
$G=\SO_0(2,n)$ with $n\geq 3$, in which cases $\pi_1(H)\cong \Z
\oplus \Z_2$). This gives an integer invariant known as the Toledo
invariant. If $\mathcal{M}_d(G)$ is the moduli space of polystable
$G$-Higgs bundles with fixed topological class of $E$, the theorem
of Corlette \cite{corlette:1988} and an adaptation of the arguments
of Simpson in \cite{simpson:1992} identify $\mathcal{R}_d(G)\cong
\mathcal{M}_d(G)$ as real analytic varieties. In a direct
generalization of the result for $G=\SL(2,\R)$, these moduli spaces
are non-empty only for values of $d$ satisfying an inequality of
Milnor--Wood type (proved by Domic and Toledo
\cite{domic-toledo:1987} and also by Turaev \cite{turaev:1984} for
the symplectic group) given by
$$
|d|\leq \rk(G/H)(g-1),
$$
where $\rk(G/H)$ is the rank of the symmetric space $G/H$. For the
semisimple classical groups defining irreducible Hermitian symmetric
spaces, namely $G=\SU(p,q), \Sp(2n,\R)$, $\SO^*(2n)$ and $\SO_0(2,n)$, we
show that this inequality is a consequence of semistability of the Higgs
bundle.

Representations with maximal Toledo invariant, so-called \emph{maximal
  representations}, are of particular geometric interest, as is
already clear in the case of $G=\SL(2,\R)$ where, as mentioned above,
they are just the uniformizing representations.  Maximal
representations have been the subject of extensive study by Burger,
Iozzi and Wienhard
\cite{burger-iozzi-wienhard:2003,wienhard:2005,burger-iozzi-wienhard:2005}
using methods of bounded cohomology. Among other things they show that
maximal representations are discrete and faithful, generalizing
Goldman's theorem for $\SL(2,\R)$ mentioned above. Another important
result proved by them is that any maximal representation is in fact
reductive, thus the restriction to reductive representations inherent
in the Higgs bundle approach is unnecessary in the case of maximal
representations.

Our main concern in this paper are maximal representations.  In this
case one finds that the geometry of $G/H$ is important. In particular,
the Shilov boundary $\Sh$ of the realization of $G/H$ as a bounded
symmetric complex domain plays a key role. There are two cases to
consider, depending on whether $G/H$ is or is not of tube-type. In the
first case, $G/H$ can be realized as tube domain over a symmetric cone
$\Omega=G'/H'$.  Moreover, the Shilov boundary $\Sh$ is symmetric
space of compact type, and the homogeneous space $G'/H'$ is the
non-compact dual of $\Sh$. This happens for the groups $G=\SU(n,n)$
$G=\Sp(2n,\R)$, $G=\SO^*(2n)$ with $n$ even and $G=\SO_0(2,n)$. In
these cases the moduli space $\mathcal{M}_d(G)$ with $d$ maximal can
be identified with another moduli space related to $G'$. This
correspondence which we call Cayley correspondence, allows us to
detect new topological invariants for the maximal representations.

In the non-tube cases (viz.\ $G=\SU(p,q)$ with $p\neq q$ and $G=\SO^*(2n)$
with $n$ odd), let $\widetilde{G}\subset G$ define the maximal
tube-type space isometrically embedded in $G/H$.  Then any maximal
representation reduces to a representation in the normalizer
$N_G(\widetilde{G})$.  This leads to a description of the moduli space
$\mathcal{M}_{d}(G)$ with $d$ maximal as a fibration, whose fibres are
isomorphic to the moduli space of $\widetilde{G}$-Higgs bundles with
maximal Toledo invariant, and whose base is the moduli space of polystable
holomorphic ${H''}^\C$-bundles, for a certain compact group
$H''$, which can de defined in terms of the Shilov boundaries of $G$ and
$\widetilde{G}$. This generalizes the rigidity
results of Toledo \cite{toledo:1989} for $p=1$ and Hern\'andez
\cite{hernandez:1991} for $p=2$ (and in the case $p>2$ for
representations satisfying a certain non-degeneracy condition).  These
two results were generalized by the authors to any reductive
representation and arbitrary $p$ in
\cite{bradlow-garcia-prada-gothen:2001,bradlow-garcia-prada-gothen:2003}
and shortly afterwards it was shown in general by Burger, Iozzi and
Wienhard \cite{burger-iozzi-wienhard:2003} that any maximal
representation stabilizes a maximal tube type subdomain of $G/H$. It
was this latter result that made us aware of the importance of the
tube type condition for the study of maximal representations.

Finally, we count the number of connected components of the moduli
spaces $\mathcal{R}_d(G)$ when $d$ corresponds to maximal values of
the Toledo invariant.  The results are summarized in
Table~\ref{tab:components} and contain as special cases results
obtained for specific groups by Goldman \cite{goldman:1988}, Hitchin
\cite{hitchin:1987a}, Xia \cite{xia:2000,xia:2003},
Xia--Markman \cite{markman-xia:2002}, and in joint work involving
Mundet i Riera and the authors
\cite{garcia-prada-mundet:2004,gothen:2001,bradlow-garcia-prada-gothen:2003,garcia-prada-gothen-mundet:2005}
(in fact many of these references give a complete count, for all
values of the Toledo invariant).

The method we use to make the final count
of components is based on the Morse function defined by Hitchin, using
the $L^2$-norm of the Higgs field. The key steps involve
characterizing the subvariety of local minima of this function and
counting the number of components of these.

\subsubsection*{Acknowledgements}

The authors thank Nigel Hitchin, Bill Goldman, Ignasi Mundet i Riera,
Domingo Toledo, Fran\c{c}ois Labourie, S. Ramanan, Joseph Wolf, Marc
Burger, Alessandra Iozzi, and Anna Wienhard for numerous useful
conversations and shared insights.  Oscar Garc\'{\i}a-Prada thanks
IHES for its hospitality and support.


\section{Surface group representations and
$G$-Higgs bundles}
\label{sec:reps-G-Higgs}

\subsection{Surface group representations}
\label{sec:representations}

Let $X$ be a closed oriented surface of genus $g$ and let
\begin{displaymath}
  \pi_{1}(X) = \{ a_{1},b_{1}, \dotsc, a_{g},b_{g} \suchthat
  \prod_{i=1}^{g}[a_{i},b_{i}] = 1 \}
\end{displaymath}
be its fundamental group.  Let $G$ be a connected reductive real  Lie group.
By a \emph{representation} of $\pi_1(X)$ in
$G$ we understand a homomorphism $\rho\colon \pi_1(X) \to G$.
The set of all such homomorphisms,
$\Hom(\pi_1(X),G)$, can be naturally identified with the subset
of $G^{2g}$ consisting of $2g$-tuples
$(A_{1},B_{1}\dotsc,A_{g},B_{g})$ satisfying the algebraic equation
$\prod_{i=1}^{g}[A_{i},B_{i}] = 1$.  This shows that
$\Hom(\pi_1(X),G)$ is a real analytic  variety, which is algebraic
if $G$ is algebraic.

The group $G$ acts on $\Hom(\pi_1(X),G)$ by conjugation:
\[
(g \cdot \rho)(\gamma) = g \rho(\gamma) g^{-1}
\]
for $g \in G$, $\rho \in \Hom(\pi_1(X),G)$ and
$\gamma\in \pi_1(X)$. If we restrict the action to the subspace
$\Hom^+(\pi_1(X),G)$ consisting of {\em reductive representations},
the orbit space is Hausdorff.  By a reductive representation we mean
one that composed with the adjoint representation in the Lie algebra
of $G$ decomposes as a sum of irreducible representations (when $G$ is
compact every representation is reductive).
Define the
\emph{moduli space of representations} of $\pi_1(X)$ in $G$
to be the orbit space
\[
\mathcal{R}(G) = \Hom^{+}(\pi_1(X),G) / G \]
with the quotient topology.

Given a representation $\rho\colon\pi_{1}(X) \to
G$, there is an associated flat $G$-bundle on
$X$, defined as
\begin{math}
  E_{\rho} = \widetilde{X}\times_{\rho}G
\end{math},
where $\widetilde{X} \to X$ is the universal cover and $\pi_{1}(X)$ acts
on $G$ via $\rho$.
This gives in fact  an identification between the set of equivalence classes
of representations  $\Hom(\pi_1(X),G) / G$ and the set of equivalence classes
of flat $G$-bundles, which in turn is parameterized by the  cohomology set
$H^1(X,G)$. We can then assign  a topological invariant  to a representation
$\rho$ given by the  characteristic class $c(\rho):=c(E_{\rho})\in \pi_1(G)$
corresponding
to $E_{\rho}$. To define this, let $\widetilde G$ be the universal covering group
of $G$. We have an exact  sequence
$$
1 \lra \pi_1(G)\lra \widetilde G \lra G \lra 1
$$
which gives rise to the (pointed sets) cohomology sequence
\begin{equation}\label{characteristic}
H^1(X, {\widetilde G}) \lra H^1(X, {G})\stackrel{c}{\lra}   H^2(X, \pi_1(G)).
\end{equation}

Since $\pi_1(G)$ is abelian, we have
$$
 H^2(X, \pi_1(G))\cong  \pi_1(G),
$$
and $c(E_\rho)$ is defined as the image of $E$ under the last map in
(\ref{characteristic}). Thus the class $c(E_\rho)$ measures  the
obstruction to lifting $E_\rho$ to a flat $\widetilde G$-bundle, and hence to
lifting $\rho$ to a representation of $\pi_1(X)$ in $\widetilde G$.  For a  fixed
 $d\in \pi_1(G)$, the
{\em moduli space of reductive representations} $\mathcal{R}_d(G)$
with topological invariant $d$
is defined as  the subvariety
$$
\mathcal{R}_d(G):=\{\rho \in \mathcal{R}(G)\; : \; c(\rho)=d\}.
$$

\subsection{$G$-Higgs bundles}
\label{sec:higgs-bundles}

Let now $X$ be a compact Riemann surface and let $K$ be its canonical
line bundle.
Let $G$ be a connected reductive real  Lie group. Let $H\subset G$ be a
maximal compact subgroup, and  $\lieg=\lieh +\liem$ be the Cartan
decomposition of $\lieg$.
A $G$-{\em Higgs bundle} over $X$ is
a pair $(E,\varphi)$ consisting of a principal holomorphic $H^\C$-bundle
$E$ over $X$ and a holomorphic section of $E(\liemc)\otimes K$, i.e.
$\varphi\in H^0(X, E(\liemc)\otimes K)$, where
$E(\liemc)$ is the
bundle  associated to $E$ via the isotropy representation of $H^\C$ in
$\liemc$.

If  $G=\GL(n,\C)$ (with its underlying real structure), we recover the
original notion of Higgs bundle
introduced by Hitchin \cite{hitchin:1987a},  consisting of
a holomorphic vector bundle $\E:=E(\C^n)$ --- associated to a principal
$\GL(n,\C)$-bundle $E$ via the standard representation ---
and a homomorphism
$$
\Phi: \E\lra \E\otimes K.
$$

The Higgs bundle $(\E,\Phi)$ is said to be {\em stable}
if
\begin{equation}\label{higgs-stability}
\frac{\deg \E'}{\rank \E'}< \frac{\deg \E}{\rank \E}
\end{equation}
for every proper subbundle $\E'\subset \E$ such that
$\Phi(\E')\subset \E'\otimes K$.
The Higgs bundle $(\E,\Phi)$ is
{\em polystable} if $(\E,\Phi)=\oplus_i (\E_i,\Phi_i)$ where
$(\E_i,\Phi_i)$ is a  stable Higgs bundles  and
${\deg \E_i}/{\rank \E_i}= {\deg \E}/{\rank \E}$.
The {\em moduli space of polystable Higgs bundles}
$\mathcal{M}(n,d)$ is defined
as the set of isomorphism classes of polystable Higgs bundles $(\E,\Phi)$
with $\rank \E=n$ and $\deg \E=d$.
Another important concept is that of {\em semistability} which
is defined by replacing the strict inequality in (\ref{higgs-stability})
by the weaker inequality.  It is immediate that
polystability implies semistability.

Similarly, there is a notion of stability, semistability and polystability for
$G$-Higgs bundles. If $G\subset \GL(n,\C)$ is a classical group,
to a $G$-Higgs bundle we can naturally associate a $\GL(n,\C)$-Higgs bundle.
The polystability of a  $G$-Higgs bundle is in fact equivalent to the
polystability of the corresponding $\GL(n,\C)$-Higgs bundle.
However, a $G$-Higgs bundle can be stable as
a $G$-Higgs bundle but not as a $\GL(n,\C)$-Higgs bundle
(\cite{bradlow-garcia-prada-gothen:2003}). Now,
 topologically, $H^\C$-bundles $E$  on $X$ are classified by
a characteristic class $d=c(E)\in \pi_1(H^\C)=\pi_1(H)=\pi_1(G)$, and
for a fixed  such  class  $d$, the
{\em moduli space of polystable $G$-Higgs bundles} $\mathcal{M}_d(G)$
is defined as  the set of isomorphism classes of polystable
$G$-Higgs bundles $(E,\varphi)$   such that $c(E)=d$.

\subsection{Correspondence of moduli spaces}
\label{sec:correspondence}

We assume now that $G$ is semisimple.
With the notation of the  previous sections, we have the following.

\begin{theorem}\label{na-Hodge}
Let $G$ be a connected non-compact semisimple real Lie group. There is a homeomorphism
$\mathcal{R}_d(G) \cong \mathcal{M}_d(G)$.
\end{theorem}

\begin{remark}
This correspondence is in fact an isomorphism of real analytic varieties.
\end{remark}

\begin{remark} There is a similar correspondence when $G$ is reductive,
replacing the fundamental group of $X$ by its universal central  extension.
\end{remark}

The proof of Theorem \ref{na-Hodge} is the combination of two
existence theorems for gauge-theoretic equations. To explain this, let
$\bE_G$ be a $C^\infty$ principal $G$-bundle over $X$ with fixed
characteristic class $d\in\pi_1(G)=\pi_1(H)$. Let $D$ be a
$G$-connection on $\bE_G$ and let $F_D$ be its curvature. If $D$ is
flat, i.e.\ $F_D=0$, then the holonomy of $D$ around a closed loop in
$X$ only depends on the homotopy class of the loop and thus defines a
representation of $\pi_1(X)$ in $G$.  This gives an
identification\footnote{even when $G$ is complex algebraic, this is
  merely a real  \emph{analytic} isomorphism, see Simpson
  \cite{simpson:1992,simpson:1994,simpson:1995}},
$$
\cR_d(G)\cong
\{\mbox{Reductive $G$-connections}\;\; D \suchthat
F_D=0\}/\GGG,
$$
where, by definition, a flat connection is reductive if the corresponding
representation of $\pi_1(X)$ in $G$ is reductive,
and $\GGG$ is the group of automorphisms of $\bE_G$  ---
the {\em gauge group}.
Let $h_0$ be a fixed  reduction
of $\bE_G$ to a $C^\infty$  $H$-bundle $\bE_H$. Every $G$-connection $D$ on $\bE_G$
decomposes uniquely as
$$
D=d_A + \psi,
$$
where $d_A$ is an $H$-connection on $\bE_H$ and
$\psi\in \Omega^1(X,\bE_H(\liem))$.  Let
$F_A$ be the curvature of $d_A$.
We consider the following set of equations for the pair $(d_A,\psi)$:
\begin{equation}\label{harmonicity}
\begin{array}{l}
F_A +\frac{1}{2}[\psi,\psi]= 0\\
d_A\psi=0  \\
d_A^\ast\psi=0.
\end{array}
\end{equation}
These equations are invariant under the action of $\HHH$, the gauge group of
$\bE_H$. A theorem of Corlette \cite{corlette:1988}, and
Donaldson \cite{donaldson:1987} for $G=\SL(2,\C)$,  says the following.
\begin{theorem}\label{corlette} There is a homeomorphism
$$
\{\mbox{Reductive $G$-connections}\;\; D \suchthat
F_D=0\}/\GGG\cong \{(d_A,\psi)\;\;\mbox{satisfying}\;\;
(\ref{harmonicity})\}/\HHH.
$$
\end{theorem}

The first two equations in (\ref{harmonicity}) are equivalent
to the flatness of $D=d_A+\psi$, and Theorem \ref{corlette}
simply says that in the $\GGG$-orbit of a reductive flat $G$-connection
$D_0$ we can find a flat $G$-connection $D=g(D_0)$ such that if we
write $D=d_A+\psi$, the
additional condition $d_A^\ast\psi=0$ is satisfied. This can be interpreted
more geometrically in terms of the reduction  $h=g(h_0)$ of $\bE_G$
to an $H$-bundle obtained by the action of $g\in\GGG$ on $h_0$.
Equation
$d_A^\ast\psi=0$ is equivalent to the harmonicity of the
$\pi_1(X)$-equivariant map $\widetilde X \ra G/H$ corresponding to
the new reduction of structure group $h$.

To establish the link  with Higgs bundles, we consider the $H$-bundle
$\bE_H$ and  the moduli space of solutions to the {\em Hitchin's equations}
for a pair $(d_A,\varphi)$ consisting
of  an $H$-connection $d_A$ and $\varphi\in \Omega^{1,0}(X,\bE_H(\liemc))$:
\begin{equation}\label{hitchin}
\begin{array}{l}
F_A -[\varphi,\tau(\varphi)]= 0\\
\dbar_A\varphi=0.
\end{array}
\end{equation}
Here $\dbar_A$ is the $(0,1)$ part of $d_A$, which defines a holomorphic
structure on $\bE_H$, and $\tau$ is the conjugation on $\liegc$ defining its
compact form. The gauge group of $\bE_H$ acts on the space of solutions
defining the moduli space of solutions. A theorem of Hitchin
\cite{hitchin:1987a} for $G=\SL(2,\C)$ and Simpson
\cite{simpson:1992} for an arbitrary semisimple complex Lie group $G$ can be  adapted to a
semisimple real Lie group $G$ \cite{bradlow-garcia-prada-mundet:2003} to give
the following.

\begin{theorem} \label{hk}
There is a homeomorphism 
$$
\cM_d(G)\cong \{ (d_A,\varphi)\;\;\mbox{satisfying}\;\;
(\ref{hitchin})\}/\HHH.
$$
\end{theorem}

To explain this correspondence we interpret the moduli
space of $G$-Higgs bundles in terms of pairs $(\dbar_E, \varphi)$ consisting
of a $\dbar$-operator on the $H^\C$-bundle $\bE_{H^\C}$ obtained from
$\bE_H$ by the extension of structure group $H\subset H^\C$, and
$\varphi\in \Omega^{1,0}(X,\bE_{H^\C}(\liemc))$
satisfying $\dbar_E\varphi=0$.
Such pairs are in correspondence with  $G$-Higgs bundles $(E,\varphi)$,
where $E$ is the holomorphic $H^\C$-bundle defined by the operator
$\dbar_E$ on $\bE_{H^\C}$ and $\dbar_E\varphi=0$ is equivalent
to $\varphi\in H^0(X,E(\liemc)\otimes K)$.
The moduli space of polystable $G$-Higgs bundles $\cM_d(G)$ can now
be identified with the orbit space
$$
\{\mbox{Polystable}\;\; (\dbar_E,\varphi)\suchthat \dbar_E\varphi=0\}/\HHH^\C,
$$
where $\HHH^\C$ is the gauge group of $\bE_{H^\C}$, which is in fact
the complexification of $\HHH$.
Since  there is a one-to-one correspondence between
$H$-connections on $\bE_H$ and $\dbar$-operators on $\bE_{H^\C}$,
the correspondence given in Theorem \ref{hk} can be interpreted
by saying that in the $\HHH^\C$-orbit of a polystable $G$-Higgs
bundle $(\dbar_{E_0},\varphi_0)$ we can find another Higgs bundle
$(\dbar_E,\varphi)$
whose corresponding pair $(d_A,\varphi)$ satisfies
$F_A -[\varphi,\tau(\varphi)]= 0$, and this is unique up to $H$-gauge
transformations.

To complete the circle, leading to Theorem \ref{na-Hodge}, we just need the
following.
\begin{proposition}\label{prop:circle}
The correspondence $(d_A,\varphi)\mapsto (d_A,\psi:=\varphi-\tau(\varphi))$
defines  a homeomorphism
$$
\{(d_A,\varphi)\;\;\mbox{satisfying}\;\;
(\ref{hitchin})\}/\HHH\cong
\{(d_A,\psi)\;\;\mbox{satisfying}\;\;
(\ref{harmonicity})\}/\HHH.
$$
\end{proposition}

For the benefit of the reader we have outlined the correspondences
between the various moduli spaces explained in this Section in a
diagram \eqref{eq:diagram-moduli}  below.
The notation used is as follows:
\begin{itemize}
\item $\bE_{H^{\C}}\to X$ is a smooth $H^{\C}$-bundle,
  \begin{itemize}
  \item $\dbar_E$ is a holomorphic structure on $\bE_{H^{\C}}$,
  \item $\phi\in \Omega^{1,0}(\bE_{H^{\C}}(\liem^{\C})$, and
  \item $\HHH^{\C}$ is the gauge group of $\bE_{H^{\C}}$;
  \end{itemize}
\item $\bE_{G}\to X$ is a smooth $G$-bundle,
  \begin{itemize}
  \item $D$ is a connection on $\bE_{G}$, and
  \item $\GGG$ is the gauge group of $\bE_{G}$;
  \end{itemize}
\item $\bE_{H}\to X$ is a smooth $H$-bundle (a reduction of
$\bE_{H^{\C}}$ and $\bE_{G}$) ,
  \begin{itemize}
  \item $A$ is a connection on $\bE_{H}$,
  \item $\psi\in \Omega^{1}(\bE_{H}(\liem)$, and
  \item $\HHH$ is the gauge group of $\bE_{H}$;
  \end{itemize}
\end{itemize}
and the maps $(1)$, $(2)$ and $(3)$ are homeomorphisms given by
\begin{itemize}
\item[(1)] $(d_A,\varphi)\mapsto (\dbar_A,\varphi)$
  \cite{hitchin:1987a,simpson:1988,bradlow-garcia-prada-mundet:2003};
\item[(2)] $(d_A,\varphi)\mapsto (d_A,\psi=\varphi-\tau(\varphi))$;
\item[(3)] $(d_A,\Psi) \mapsto D = d_A + \psi$
  \cite{hitchin:1987a,corlette:1988}.
\end{itemize}

The following diagram outlines the correspondences between the various
moduli spaces explained in this Section (cf.\ Theorems~\ref{na-Hodge},
\ref{corlette} and \ref{hk} and Proposition~\ref{prop:circle}).

\newlength{\wi}
\newlength{\wii}
\newlength{\wiii}
\newlength{\wiiii}
\settowidth{\wi}{polystable}
\settowidth{\wii}{$F_A-[\varphi,\tau(\varphi)]=0$}
\settowidth{\wiii}{reductive}
\settowidth{\wiiii}{$F_A +\tfrac{1}{2}[\psi,\psi]=0$}

\begin{equation}\label{eq:diagram-moduli}
  \begin{CD}
    \mathcal{M}_d(G)
    \cong
    \left\{ (\dbar_E,\varphi)\suchthat
    \parbox{\wi}{\centering $\dbar_E\varphi=0$ \\
      polystable} \;\right\} {/ \HHH^{\C}}
      @<<(1)<
   \left\{ (d_A,\varphi)\suchthat
    \parbox{\wii}{\centering
      $\begin{aligned}
        \dbar_A\varphi &=0 \\
        F_A-[\varphi,\tau(\varphi)] &= 0
      \end{aligned}$}
  \;\right\} {/ \HHH} \\
  {\rotatebox{90}{$\cong$}}@. @VV(2)V \\
 \mathcal{R}_d(G) \cong \left\{ D \suchthat
    \parbox{\wiii}{\centering $F_D =0$
      \\ reductive }\;\right\} {/ \GGG}
    @<<(3)<
  \left\{ (d_A,\psi)\suchthat
    \parbox{\wiiii}{\centering $
      \begin{aligned}
        F_A +\tfrac{1}{2}[\psi,\psi] &= 0\\
        d_A\psi &=0  \\
        d_A^\ast\psi &=0
      \end{aligned}$}
  \;\right\} {/ \HHH}.
 \end{CD}
\end{equation}

\section{Isometry groups of Hermitian symmetric spaces}
\label{sec:isom-groups-herm}

\subsection{$G$-Higgs bundles for the classical Hermitian
symmetric spaces}
\label{sec:G-higgs-bundles}

Our goal in this paper is to study representations
of the fundamental group in the case in which $G/H$ is a Hermitian
symmetric space.
This means that  $G/H$ admits a complex structure compatible with
the Riemannian structure of $G/H$, making $G/H$  a K\"ahler manifold.
If $G/H$ is irreducible, the centre of $\lieh$ is one-dimensional
and the almost complex structure
on $G/H$ is defined by a generating element in $J\in Z(\lie{h})$
(acting through the isotropy representation on $\lie{m}^{\C}$). This
complex structure defines a decomposition

\begin{equation}\nonumber
\lie{m}^\C=\lie{m}_++\lie{m}_-,
\end{equation}
where $\liem_+$ and $\liem_-$ are the $(1,0)$ and the $(0,1)$ part of
$\liem^\C$ respectively. Table~\ref{tab:HSS} (see
Sec.~\ref{sec:tables}) shows the main ingredients
for the irreducible classical Hermitian symmetric spaces.

Let now $(E,\varphi)$ be a $G$-Higgs bundle over a compact Riemann
surface $X$. The decomposition $\liemc=\liemp+\liemm$
gives a vector bundle decomposition $E(\liemc)= E(\liemp) \oplus E(\liemm)$
and hence
$$
\varphi=(\beta, \gamma)\in H^0(X,E(\liemp)\otimes K)\oplus
H^0(X,E(\liemm)\otimes K)= H^0(X,E(\liemc)\otimes K).
$$

In Table~\ref{tab:higgs-HSS}, we describe the $G$-Higgs bundles
for the various groups appearing in Table~\ref{tab:HSS}. It is
sometimes convenient
to replace the $H^\C$-bundle $E$ by a vector bundle associated to the
standard representation of $H^\C$. For example, for $G=\SU(p,q)$,
$$
H^\C=S(\GL(p,\C)\times \GL(q,\C))
:=\{(A,B)\in  \GL(p,\C)\times \GL(q,\C)\suchthat \det B=(\det A)^{-1}\},
$$
and the $H^\C$-bundle $E$ is replaced by two holomorphic vector
bundles $V$ and $W$ or rank $p$ and $q$, respectively such that
$\det W= (\det V)^{-1}$.

As mentioned in Sec.~\ref{sec:higgs-bundles}, via the natural inclusion
$G\subset G^\C\subset \SL(N,\C)$ for $G$ in Table~\ref{tab:higgs-HSS},
to a $G$-Higgs bundle $(E,\varphi)$ we can naturally associate an
$\SL(N,\C)$-Higgs bundle $(\E,\Phi)$, where $\E$ is a holomorphic vector
bundle with trivial determinant. This is a very useful correspondence
that we also  describe in Table~\ref{tab:higgs-HSS}.

\subsection{Toledo invariant and Milnor--Wood inequalities}
\label{sec:toledo-MW}

Let $G$ be a semisimple Lie group such that  $G/H$ is an irreducible
Hermitian symmetric space. Then the torsion-free  part of $\pi_1(H)$
is isomorphic to $\Z$ and hence the topological invariant of either
a representation of $\pi_1(X)$ in $G$, or of a $G$-Higgs bundle, is
measured by an integer $d\in\Z$, known as the {\em Toledo
invariant}. In the classical cases described in Table~\ref{tab:HSS}
this coincides with the degree of a certain vector bundle. In fact,
besides $G=\SO_0(2,n)$ with $n\geq 3$, for which $\pi_1(H)\cong \Z
\oplus \Z_2$, for all the other groups in Table~\ref{tab:HSS}
$\pi_1(H)\cong \Z$.

Let $G=\SU(p,q)$. As we can see from Table~\ref{tab:higgs-HSS},
an $\SU(p,q)$-Higgs
bundle over $X$ is defined by  a 4-tuple $(V,W,\beta,\gamma)$ consisting
of holomorphic vector bundles $V$ and $W$ of rank $p$ and $q$, respectively,
such that $\det W={(\det V)}^{-1}$, and
homomorphisms
$$
\beta: W\lra V \otimes K \;\;\;\mbox{and}\;\;\;
\gamma: V\lra W \otimes K.
$$
If $(V,W,\beta,\gamma)$ is polystable then the associated Higgs
bundle $(\E,\Phi)$ with
$$
\E=V\oplus W\;\;\; \mbox{and}\;\;\;
\Phi =
\begin{pmatrix}
  0 & \beta \\
  \gamma  & 0
\end{pmatrix}
$$
is also  polystable and in particular  semistable.
Applying the semistability numerical criterion to special Higgs subbundles
defined by the kernel and image of $\Phi$ (see
\cite{bradlow-garcia-prada-gothen:2003}) we obtain
\begin{align}
d&\leq \rank(\gamma) (g-1) \label{semistability}\\
-d&\leq \rank(\beta) (g-1),
\end{align}
which gives the inequality
\begin{equation}\label{milnor-wood}
 |d|\leq \min\{p,q\}(g-1),
\end{equation}
that  generalizes the inequality of Milnor--Wood for $G=\SU(1,1)$.
Using similar arguments  one can show the various {\em Milnor--Wood type
inequalities} for  $G$ in Table~\ref{tab:higgs-HSS}.
In fact, the
inequalities for $G=\Sp(2n,\R)$ and $G=\SO^*(2n)$ can be obtained
from that of $\SU(n,n)$, via the natural inclusion of these two
groups in $\SU(n,n)$.
These inequalities for representations
of $\pi_1(X)$ in $G$ have been proved by Domic and Toledo \cite{domic-toledo:1987}
using other methods.  One can observe that the bound for the Toledo
invariant can uniformly written as $\rank(G/H)(g-1)$.

\begin{remark}\label{rem:positive-toledo}
  Duality gives an isomorphism $\cM_{d}(G)\cong \cM_{-d}(G)$ for every
$G$ in Table~\ref{tab:higgs-HSS}. For example, for $\SU(p,q)$ this
isomorphism is defined by the map
$$
(V,W,\beta,\gamma)\mapsto (V^*,W^*,\gamma^t,\beta^t).
$$
There is hence  no loss of generality in considering only the case with positive
Toledo invariant, which we will do from now on.
\end{remark}

Our main interest in this paper is the case when the Toledo invariant $d$
is maximal, that is $|d|=d_{\max}$. Let
$$
\cM_{\max}(G):=\cM_{d_{\max}}(G)
$$
be the moduli space of $G$-Higgs bundles with maximal Toledo invariant,
and similarly $\cR_{\max}(G):=\cR_{d_{\max}}(G)$ for the corresponding
moduli space of representations.

Our goal in the following sections is to study the geometry of
maximal Higgs bundles and to count the number of
connected components of $\cM_{\max}(G)$ for $G$ in Table~\ref{tab:higgs-HSS}.
The case of $G=\SU(1,1)\cong \SL(2,\R)$ and  $G=SO_0(2,1)\cong \PSL(2,\R)$
was studied by Goldman \cite{goldman:1980}, who showed  that the
moduli space of maximal representations in
$\SL(2,\R)$ has $2^{2g}$ connected components isomorphic to
Teichm\"uller space, all of which get identified when we consider
representations in  $\PSL(2,\R)$.

\subsection{Low rank phenomena}
  \label{sec:sp-so-coincidence}
  There are certain low rank coincidences between the Higgs bundles
  given in Table~\ref{tab:higgs-HSS}, coming from the special low
  dimensional isomorphisms between Lie groups.

  The most basic instance of this is isomorphism $\SU(1,1) \cong
  \Sp(2,\R)$ which clearly gives rise to equivalent Higgs vector
  bundle data $(\E,\Phi)$, namely a line bundle $M$ (which equals the
  $V$ of Table~\ref{tab:higgs-HSS}) and a pair of sections $\beta\in
  H^0(M^2K)$ and $\gamma\in H^0(M^{-2}K)$ --- we note that this is
  also the data corresponding to a $\SL(2,\R)$-Higgs bundle (cf.\
  Hitchin \cite{hitchin:1987a}), as is to be expected from the
  isomorphism $\Sp(2,\R) \cong \SL(2,\R)$.

  In a similar way, the covering $\SU(1,1) \cong \Sp(2,\R) \cong
  \Spin_0(2,1) \to \SO_0(2,1)$ shows that a Higgs bundle for one of
  the former groups gives rise to one for the latter.  Explicitly, if
  $(M,\beta,\gamma)$ is as above, then the associated
  $\SO_0(2,1)$-Higgs bundle $(\E= V \oplus W,\Phi)$ has $V = M^2
  \oplus M^{-2}$ and $W = \cO$; thus the $L$ of
  Table~\ref{tab:higgs-HSS} is $L = M^2$.  Note that the Toledo
  invariants are related by
  \begin{displaymath}
    d_{\Spin_0(2,1)} = \deg(M) = \frac{1}{2}\deg(L) =
    \frac{1}{2}d_{\SO_0(2,1)},
  \end{displaymath}
  and that a $\SO_0(2,1)$-Higgs bundle lifts to a $\Spin_0(2,1)$-Higgs
  bundle if and only it has even Toledo invariant.

  Analogous phenomena occur for various other (local) isomorphisms,
  here we shall describe explicitly just one more such situation of particular
  interest, corresponding to the covering
  $\Sp(4,\R) \cong \Spin_0(2,3) \to \SO_0(2,3)$.
  Let $(V,\beta,\gamma)$ be the vector bundle data corresponding to a
  $\Sp(4,\R)$-Higgs bundle as in Table~\ref{tab:higgs-HSS} and let $L
  = \Lambda^2 V$ be the determinant bundle of the rank $2$ bundle
  $V$.  Then $S^2V$ has a non-degenerate quadratic form with values in
  $L^2$ defined by
  \begin{displaymath}
    Q(x\otimes y,x'\otimes y') = (x \wedge x')\otimes(y\wedge y')
  \end{displaymath}
  and thus we obtain an orthogonal bundle $(W,Q_W)$ letting
  \begin{displaymath}
    W = S^2V \otimes L^{-1}
  \end{displaymath}
  with the induced quadratic form $Q_W$. Thus we have the required
  bundle data $(W,Q_W)$ and $L$ to define a $\SO_0(2,3)$-Higgs bundle.
  Since $S^2V = W \otimes L \cong W^*\otimes L$, the section $\beta\in
  H^0(S^2V \otimes K)$ can be viewed as a section of $\Hom(W,L)\otimes
  K$ and, similarly, $\gamma$ can be viewed as a section of
  $\Hom(W,L^{-1})\otimes K$.  Note that, since $\deg(V) =
  \deg(\Lambda^2V) = \deg(L)$, in this case the Toledo invariants are
  the same.  However, the topological classification of $SO_0(2,3)$
  bundle involves \emph{two} classes, namely the degree of $L$ (i.e.\
  the Toledo invariant) \emph{and} the second Stiefel--Whitney class
  $w_2(W,Q_w) \in \Z/2$
  of the $\SO(3,\C)$ bundle $(W,Q_W)$.  One can show without too much
  difficulty that an $\SO_0(2,3)$-Higgs bundle given by the data
  $(L,W,Q_W,\beta,\gamma)$ lifts to a
  $\Sp(4,\R)$-Higgs bundle if and only if
  \begin{displaymath}
    \deg(L) = w_2(W,Q_W).
  \end{displaymath}

\begin{remark}\label{rem:so22}
  The group $\SO_0(2,2)$ is special because the associated Hermitian
  symmetric space is not irreducible; in fact $\SO_0(2,2)$ is
  isogenous to $\SL(2,\R) \times \SL(2,\R) \cong \Spin_0(2,2)$.  Of
  course, the results for irreducible Hermitian symmetric spaces given
  here, can be applied to obtain results for isometry groups of all
  Hermitian symmetric spaces. For this reason we shall exclude this
  group from our considerations in this paper.
\end{remark}

\section{Maximal Toledo invariant}
\label{sec:maximal-toledo}

\subsection{Tube type condition}
\label{sec:tube-type}
We refer to
\cite{faraut-koranyi:1994,helgason:1978,koranyi-wolf:1965,satake:1980}
for details regarding this section.

It is well-known that a Hermitian symmetric space of non-compact type
$G/H$ can be realized as a bounded symmetric domain.
For the classical groups this
is due to Cartan, while the  general case
is given by the Harish-Chandra embedding $ G/H \ra \liemp$
which defines
a biholomorphism between $G/H$ and
the bounded symmetric domain $\cD$ given  by the image of $G/H$ in
the complex  vector space $\liemp$.
Now, for any bounded domain $\cD$ there is  the {\em Shilov boundary}
of $\cD$ which is defined as the smallest closed subset $\Sh$
of the topological boundary $\partial \cD$ for which  every function
 $f$ continuous on $\overline{\cD}$ and holomorphic on $\cD$ satisfies
that
$$
|f(z)|\leq \max_{w \in \Sh} |f(w)|\;\;\mbox{for every}\;\; z\in\cD.
$$
The Shilov boundary $\Sh$ is the unique closed $G$-orbit in
$\partial \cD$.

The simplest situation to consider is that of the hyperbolic plane.
The Poincar\'e disc is its realization as a bounded symmetric domain.
However, we know that the hyperbolic plane  can also be realized as
the upper-half plane. There are other Hermitian symmetric spaces that,
like the hyperbolic plane, admit a realization similar to   the
upper-half plane. These are the tube type symmetric spaces.

Let $V$ be a real vector space and let $\Omega \subset V$ be an open cone
in $V$. A tube over the cone $\Omega$  is a domain of the form
$$
T_\Omega=\{u+iv\in V^\C, u\in V, v\in \Omega\}.
$$
A domain $\cD$ is said to be of {\em tube type}  if it is biholomorphic
to a tube $T_\Omega$. In the case of a symmetric domain the cone $\Omega$
is also symmetric.
An important characterization of the tube type symmetric domains is given
by the following.

\begin{proposition}Let $\cD$ be a bounded symmetric domain.
The following are equivalent:

(i) $\cD$ is of tube type.

(ii) $\dim_\R \Sh=\dim_\C\cD$.

(iii) $\Sh$ is a symmetric space of compact type.

\end{proposition}

There is a generalization of the Cayley map that sends the  unit disc
biholomorphically to the upper-half plane. Let  $\cD$ be  the bounded domain
associated to a Hermitian  symmetric space $G/H$.
Acting by  a particular element in $G$, known as the {\em Cayley element},
one obtains a map
$$
c: \cD\lra \liemp
$$
which is called  the  {\em Cayley transform}.
A  relevant fact for us is the following.

\begin{proposition}
Let $\cD$ be the  symmetric domain corresponding to the Hermitian
symmetric space $G/H$. Let $\cD$ be of tube type. Then the image by
the Cayley transform $c(\cD)$  is biholomorphic to a tube domain
$T_\Omega$ where the symmetric cone $\Omega$ is the non-compact
dual of the  Shilov boundary of $\cD$. In fact the Shilov boundary
is a symmetric space isomorphic to  $H/H'$ for a certain subgroup
$H'\subset H$, and $\Omega=G'/H'$ is its non-compact dual symmetric
space.
\end{proposition}

\begin{proposition}
  \begin{enumerate}
  \item The symmetric spaces defined by $\Sp(2n,\R)$, $\SO_0(2,n)$ are
    of tube type.
  \item The symmetric space defined by $\SU(p,q)$ is of tube type if
    and only if $p=q$.
  \item The symmetric space defined by $\SO^*(2n)$ is of tube type if
    and only if $n$ is even.
  \end{enumerate}
\end{proposition}

For a  tube type classical irreducible symmetric space $G/H$,
Table~\ref{tab:tube} indicates  the Shilov boundary
$\Sh=H/H'$, its non-compact dual $\Omega= G'/H'$,
the isotropy representation space $\liem'$
and its complexification ${\liem'}^\C$, corresponding to the Cartan
decomposition of the Lie algebra $\lieg'=\lieh'+\liem'$ of $G'$.
The vector space $\liem'$  has the structure of a Euclidean Jordan
algebra,  where the cone $\Omega$ is realized.



The study of certain problems in non-tube type domains
can be reduced to the tube type thanks to the following.

\begin{proposition}
Let $G/H$ be  a Hermitian symmetric space of non-compact type.
There exists a subgroup  $\widetilde G\subset G$ such that
$\widetilde G/\widetilde H\subset
G/H$ is a maximal isometrically embedded  symmetric space of tube type,
where  $\widetilde H\subset \widetilde G$ is a maximal compact subgroup.
\end{proposition}

Table~\ref{tab:non-tube} gives the maximal symmetric space of tube
type isometrically embedded in the two series of irreducible classical
symmetric spaces of non-tube type. We describe also the Shilov
boundaries of $G/H$ and $\widetilde G/\widetilde H$ which are of the form
$\Sh=H/H'$, and $\widetilde{\Sh}=\widetilde H/\widetilde H'$, respectively. Notice
that in the non-tube case the Shilov boundary $\Sh$ is a homogeneous
space $H/H'$ but it is not symmetric.



\subsection{Tube type domains and Cayley correspondence}
\label{sec:cayley}
It turns out that the behaviour of maximal
representations and $G$-Higgs bundles is governed by the tube type
nature of $G/H$ and the  geometry of its
Shilov boundary. To explain this, let $L$ be a holomorphic
 line bundle over $X$.
Let $G$ be a real reductive  Lie group. Let $H\subset G$ be a
maximal compact subgroup, and  $\lieg=\lieh +\liem$ be the Cartan
decomposition of $\lieg$.
An  {\em $L$-twisted $G$-Higgs pair} is
a pair $(E,\varphi)$ consisting of a principal holomorphic $H^\C$-bundle
$E$ and a holomorphic section of $E(\liemc)\otimes L$, where
$E(\liemc)$ is the bundle  associated to $E$ via the isotropy  representation
of $H^\C$ in $\liem^\C$. Note that a $G$-Higgs bundle is simply a
$K$-twisted $G$-Higgs pair. Let $\mathcal{M}_L(G)$ be the moduli space
of polystable $L$-twisted $G$-Higgs pairs.

\begin{theorem} \label{cayley-correspondence}
Let $G$ be a connected semisimple classical Lie group such that $G/H$ is a
Hermitian symmetric space of tube type,
and let $\Omega= G'/H'$ be the non-compact dual of the Shilov boundary
$\Sh=H/H'$ of $G/H$.
Then
\begin{equation}\label{cayley}
\mathcal{M}_{\max} (G) \cong  \mathcal{M}_{K^2}(G').
\end{equation}
\end{theorem}
In analogy with the Cayley transform of the previous section, we
call the isomorphism
given in Theorem \ref{cayley-correspondence} {\em Cayley correspondence}.
To prove this correspondence, it suffices to do it for all $G$ in Table
\ref{tab:tube}. We will sketch the main arguments case by case.
See \cite{gothen:2001, bradlow-garcia-prada-gothen:2001,
bradlow-garcia-prada-gothen:2003,
garcia-prada-mundet:2004,
garcia-prada-gothen-mundet:2005,
bradlow-garcia-prada-gothen:preparation-2}
for details.

\begin{itemize}

\item $G=\SU(n,n)$:

\end{itemize}

An  $\SU(n,n)$-Higgs bundle over $X$ is defined by  a 4-tuple
$(V,W,\beta,\gamma)$ consisting
of two  holomorphic vector bundles $V$ and $W$ of rank $n$
such that $\det W={(\det V)}^{-1}$, and homomorphisms
$$
\beta: W\lra V \otimes K \;\;\;\mbox{and}\;\;\;
\gamma: V\lra W \otimes K.
$$

Suppose that the Toledo invariant $d=\deg V$ is maximal and positive,
that is,  $d=n(g-1)$. From (\ref{semistability}) we deduce that
$\gamma$ must be an isomorphism. Let
$\theta:  W\ra W \otimes K^2$  be defined as
$\theta= (\gamma\otimes I_K) \circ \beta $, where $I_K:K\ra K$ is the
identity map.

The condition $\det W={(\det V)}^{-1}$, together with the
isomorphism $\gamma$ imply that $(\det W)^2\cong K^{-n}$.
Now, if we choose a square root of the canonical bundle,
$L_0=K^{1/2}$, and define $\widetilde W=W\otimes L_0$, we have that
${(\det \widetilde W)}^2=\cO$ and hence the structure group of $\widetilde W$ is
the kernel of the group homomorphism $\GL(n,\C)\ra \C^*$ given by
$A\mapsto {(\det A)}^2$. This kernel is isomorphic to the semidirect
product $\SL(n,\C) \rtimes \Z_2$, where $\Z_2=\{\pm I\}$ and has then
two connected components.
The choice of a 2-torsion element in the
Jacobian of $X$ for  $\det \widetilde W$  defines an invariant
that takes $2^{2g}$ values.

Let $\widetilde\theta: \widetilde W\ra \widetilde W \otimes K^2$  be defined as
$\widetilde\theta=\theta\otimes I_{L_0}$. The map
\begin{equation}
  \label{eq:cayley-su}
  (V,W,\beta,\gamma)\mapsto (\widetilde W, \widetilde \theta)
\end{equation}
gives the bijection
(\ref{cayley}) with $G'= \SL(n,\C) \rtimes \Z_2$ (see Table~\ref{tab:tube}).

\begin{itemize}

\item  $G=\Sp(2n,\R)$:

\end{itemize}

A $\Sp(2n,\R)$-Higgs bundle
over $X$ is defined by  a triple
$(V,\beta,\gamma)$ consisting
of a rank $n$  holomorphic vector bundles $V$
and {\em symmetric}  homomorphisms
$$
\beta: V^*\lra V \otimes K \;\;\;\mbox{and}\;\;\;
\gamma: V\lra V^* \otimes K.
$$

If the Toledo invariant $d=\deg V$ is maximal and positive,
that is,  $d=n(g-1)$,  again from  (\ref{semistability}) we deduce that
$\gamma$ is   an isomorphism. Let $L_0=K^{1/2}$ be a fixed  square
root of $K$, and define $ W=V^*\otimes L_0$. Then
$Q:=\gamma\otimes I_{L_0^{-1}}: W^* \ra W$ is a symmetric isomorphism
defining an orthogonal structure on $W$, in other words,  $(W,Q)$ is an
$\OO(n,\C)$-holomorphic bundle. The $K^2$-twisted endomorphism
$\theta:W\ra W \otimes K^2$
defined by
$\theta = (\gamma\otimes I_{K\otimes L_0})\circ \beta \otimes I_{L_0}$
is $Q$-symmetric and hence $(W,\theta)$ defines a
$K^2$-twisted $\GL(n,\R)$-Higgs pair, from which we can recover the original
$\Sp(2n,\R)$-Higgs bundle, giving  the bijection (\ref{cayley}) in this case.

\begin{itemize}

\item $G=\SO^*(2n)$, with $n=2m$:

\end{itemize}

A $\SO^*(2n)$-Higgs bundle is
over $X$ is defined by  a triple
$(V,\beta,\gamma)$ consisting
of a rank $n$  holomorphic vector bundles $V$
and skew-symmetric  homomorphisms
$$
\beta: V^*\lra V \otimes K \;\;\; \mbox{and}\;\;\;
\gamma: V\lra V^* \otimes K.
$$

Since $n=2m$ is even, the maximal value of the Toledo invariant
(see Table~\ref{tab:higgs-HSS}) is $d_{\max}=n(g-1)$. If $d=n(g-1)$ then,
as in the
previous cases, $\gamma$ is  an isomorphism, and if  $L_0=K^{1/2}$ is
 a fixed  square
root of $K$, and we define $ W=V^*\otimes L_0$, the homomorphism
$\omega:=\gamma\otimes I_{L_0^{-1}}: W^* \lra W$ is a skew-symmetric
isomorphism
defining a symplectic  structure on $W$, that is,  $(W,\omega)$ is a
$\Sp(2m,\C)$-holomorphic bundle. The $K^2$-twisted endomorphism
$\theta:W\lra W \otimes K^2$
defined by
$\theta = (\gamma\otimes I_{K\otimes L_0})\circ \beta \otimes I_{L_0}$
is in this case skew-symmetric with respect to  $\omega$ and hence
 $(W,\theta)$ defines a
$K^2$-twisted $G'$-Higgs pair for $G'=\U^*(2m)$.
The map $(V,\beta,\gamma)\mapsto (W,\theta)$
gives  the bijection (\ref{cayley}) in this case.

\begin{itemize}

\item $G=\SO_0(2,n)$:

\end{itemize}

A $\SO_0(2,n)$-Higgs bundle  is defined by  a $\SO(2,\C)$-bundle
$$
\big( V=L\oplus L^{-1},Q_V=
\begin{pmatrix}
  0 & \beta \\
  \gamma  & 0
\end{pmatrix}\big)
$$
where $L$ is a holomorphic line bundle,  and a $\SO(n,\C)$-bundle
$(W,Q_W)$, together with  homomorphisms
$$
\beta: W\lra L\otimes K\;\;\;\mbox{and}\;\;\;
\gamma: W\lra L^{-1}\otimes K.
$$
The maximal case corresponds to  $d=\deg L=2g-2$. In this situation one
can show that $\gamma$ has (maximal) rank one at all points and hence it is
surjective. If we define  $F:=\ker\gamma$, we  have a sequence
\begin{equation}\label{extension}
0 \lra F \lra W  \lra  L^{-1}\otimes K  \lra 0.
\end{equation}
One can show that this sequence splits and $F$ inherits a
$\OO(n-1,\C)$-structure.
Let $L_0:=L^{-1}\otimes K$. From the exact sequence  we deduce that
$L_0\otimes\det F\cong\cO$ and hence $L_0^2=\cO$. In other words,
$L^2\cong K^2$.
Now, according to the decomposition $W\cong F\oplus L_0$, we can decompose
$\beta=\beta'+\beta''$ with $\beta': F\lra L\otimes K$ and
$\beta'': L_0\lra L\otimes K$. Tensoring these homomorphisms by
$L_0$, we obtain $\theta': F\otimes L_0\ra K^2$ and
$\theta'': \cO\ra K^2$. The map
$$
(L,W,Q_W,\beta,\gamma)\mapsto (F,\theta',\theta'')
$$
defines the correspondence (\ref{cayley}) now.

In all the cases above we have  to show of course that
the  corresponding (poly)stability conditions in
$\cM_{\max}(G)$ and $\cM_{K^2}(G')$ are equivalent.

The correspondence (\ref{cayley}) brings to the surface
new topological  invariants of a maximal representation --- the invariants
of the $H'^\C$-bundle --- which are not {\em a priori} ``visible''.
The new invariants  will account to a certain extent for the abundance of
connected  components in most cases.

\begin{remark}
We believe that Theorem \ref{cayley-correspondence} is  true  also in the
non-classical case.
To show this it would suffice to check the case of the rank 3
irreducible  exceptional domain, which is  obtained from a real form of
$E_7$. It would be very interesting, however,  to find a proof independent
of classification theory.
\end{remark}

\subsection{Non-tube type domains and rigidity of  maximal representations}
\label{sec:non-tube-type}

We study now maximal $G$-Higgs bundles and representations when $G/H$ is
not of tube type (see Table~\ref{tab:non-tube}).
Let us start with  $G=\SU(p,q)$ with $p\neq q$.
Without loss of generality we assume that $p<q$.
 The maximal value of the Toledo
invariant is then $d_{\max}=p(g-1)$. As shown in
\cite{gothen:2001,bradlow-garcia-prada-gothen:2003}, it turns out that
there are no stable $\SU(p,q)$-Higgs bundles. In fact, every polystable
$\SU(p,q)$-Higgs bundle  $(V,W,\beta,\gamma)$ is strictly semistable
and decomposes  as a direct sum
\begin{equation}
\label{eq:rigidsupqhiggs}
(V,W,\beta,\gamma) \cong  (V,W',\beta,\gamma) \oplus  (0,W'',0,0),
\end{equation}
of a maximal polystable  $\U(p,p)$-Higgs bundle and a polystable
$\GL(q-p,\C)$-bundle with zero Chern class, 
where $\gamma\colon V \xrightarrow{\cong} W' \otimes K$, with
$W'=\im\gamma\otimes K^{-1}$ and 
$W''=W/W'$. Since 
\begin{displaymath}
  \det(V)\otimes \det(W') \otimes(W'') \cong \mathcal{O}, 
\end{displaymath}
this means that  the $\SU(p,q)$-Higgs bundle reduces to an
$\SSS(\U(p,p) \times \U(q-p))$-Higgs bundle.

We have the exact sequence
\begin{align*}
    1 \to \SU(p,p) \to \SSS(\U(p,p) \times \U(q-p)) &\to
    \SSS(\U(1)\times\U(q-p)) \to 1 \\
    (A,B) &\mapsto (\det(A),B),
\end{align*}
from which we conclude the following.

\begin{theorem}\label{rigidity-supq}
 Let $p<q$. Then the moduli space $\mathcal{M}_{\max} (\SU(p,q))$
 fibres over $M(\GL(q-p,\C))$ with
fibre isomorphic to $\mathcal{M}_{\max}(\SU(p,p))$, 
where $M(\GL(q-p,\C))$ is the moduli space of polystable vector bundles of rank
$q-p$ and zero Chern class.
\end{theorem}

Similarly, for $G=\SO^*(2n)$ with $n=2m+1$, the maximal value of
the Toledo invariant is $d_{\max}=(n-1)(g-1)$, and if $d=d_{\max}$ there are no stable
$\SO^*(4m+2)$-Higgs bundles
and every polystable $\SO^*(4m+2)$-Higgs bundle $(V,\beta,\gamma)$  decomposes as
$$
(V,\beta,\gamma)=(V',\beta,\gamma) \oplus (L,0,0),
$$
where $(V',\beta,\gamma)$ is a maximal polystable $\SO^*(4m)$-Higgs bundle,
with $V'=(\im \gamma)^* \otimes K$ and $L=V/V'$ is a line bundle of $0$ degree
(see \cite{bradlow-garcia-prada-gothen:preparation-2} for details).

We thus have the following.

\begin{theorem}\label{rigidityso*2n}
$$
\mathcal{M}_{\max} (\SO^*(4m+2))\cong \mathcal{M}_{\max}(\SO^*(4m))\times J(X),
$$
where $J(X)$ is the Jacobian of $X$.
\end{theorem}

Since, as we know from Table~\ref{tab:non-tube}, the two cases
discussed above are the only ones defining classical irreducible
Hermitian symmetric spaces of non-tube type, we conclude the
following.

\begin{theorem} \label{rigidity-higgs} Let $G$ be a connected
  semisimple classical Lie group such that $G/H$ is Hermitian
  symmetric space of non-compact type.  Let $\widetilde G\subset G$ be
  a subgroup defining the largest isometrically embedded subspace
  $\widetilde G/\widetilde H\subset G/H$, and let $H''=H'/\widetilde
  H'$ where $\Sh=H/H'$ and $\widetilde{\Sh}=\widetilde H/\widetilde
  H'$ are the Shilov boundaries of $G/H$ and $\widetilde G/\widetilde
  H$, respectively.  Then the following holds:

(1) Every $G$-Higgs bundle in
$\mathcal{M}_{\max} (G)$ is strictly polystable and reduces to a 
$N_G(\widetilde{G})$-Higgs bundle, where the
normalizer of $\tilde{G}$ in $G$, $N_G(\tilde{G})$,  fits in 
the exact sequence
\begin{displaymath}
  1 \to \tilde{G} \to N_G(\tilde{G}) \to H'' \to 1.
\end{displaymath}

(2) The moduli space $\mathcal{M}_{\max} (G)$ fibres over
$M({H''}^\C)$, with fibre  $\mathcal{M}_{\max}(\widetilde{G})$,
where $M({H''}^\C)$ is the moduli space of polystable holomorphic
${H''}^\C$-bundles with zero characteristic class.
\end{theorem}

From Theorem \ref{na-Hodge} and the theorems of Narasimhan--Seshadri
\cite{narasimhan-seshadri:1965} and Ramanathan \cite{ramanathan:1975},
which identify the moduli space of polystable holomorphic
${H''}^\C$-bundles with trivial characteristic class with the moduli
space of representations of $\pi_1(X)$ in $H''$, we obtain the
following.

\begin{theorem}\label{rigidity-reps}
  With the same hypotheses and notation as in 
{\rm Theorem \ref{rigidity-higgs}} we have the following:

(1) Every representation in $\mathcal{R}_{\max} (G)$ is 
reducible and factors through  a representation in $N_G(\widetilde{G})$.

(2) The moduli space $\mathcal{R}_{\max} (G)$ fibres over
$\mathcal{R}(H'')$, with fibre  $\mathcal{R}_{\max}(\widetilde{G})$,
where $\mathcal{R}(H'')$ is the moduli space of representations of
$\pi_1(X)$ in $H''$.

\end{theorem}

Theorem \ref{rigidity-reps} had been proved by Toledo
\cite{toledo:1989} for $G=\SU(1,q)$.  For general $p$ the result had
been proved by Hern\'andez \cite{hernandez:1991}, under a certain
non-degeneracy condition on the representation, which he was able to
show is always satisfied for $p=2$.  The result was then proved in
\cite{bradlow-garcia-prada-gothen:2001,bradlow-garcia-prada-gothen:2003}
for any reductive representation in $\SU(p,q)$.  Finally Burger, Iozzi
and Wienhard \cite{burger-iozzi-wienhard:2003} showed that any maximal
representation stabilizes a maximal tube type subdomain of
$\widetilde{G}/\widetilde{H} \subset G/H$. From this, and  their results in
\cite{burger-iozzi-wienhard:2005},   Theorem~\ref{rigidity-reps}
should follow directly for general $G$.  
On the other hand, to generalize Theorems \ref{rigidity-higgs} and 
\ref{rigidity-reps} from the Higgs bundle point of view, it suffices
to prove corresponding
results for the only non-tube rank 2 irreducible exceptional domain
which is obtained from a real form of $E_6$.

\begin{remark}
  Theorems \ref{rigidity-higgs} and \ref{rigidity-reps} establish a
  certain kind of {\em rigidity} for maximal $G$-Higgs bundles, and
  hence for surface group representations in $G$, when $G/H$ is not of
  tube type.  Namely, since the expected complex dimension of
  $\mathcal{M}_d(G)=\dim_\R \lieg (g-1)$ (as can be computed using the
  deformation theory in Sec.~\ref{sec:deformation-theory}), in the non
  tube situation $\dim_\R \widetilde{\lieg} + \dim_\R \lieh'' < \dim
  \lieg$ and hence the dimension of $\mathcal{M}_{\max}$ is smaller
  than expected (here $\lieh''$ and $\widetilde \lieg$ are the Lie
  algebras of $\widetilde G$ and $\widetilde H''$ in Theorem
  \ref{rigidity-higgs}).
\end{remark}

The fact that $M({H''}^\C)$ is connected
(\cite{narasimhan-seshadri:1965,ramanathan:1975}) leads to the
following Corollary to Theorem~ \ref{rigidity-higgs}.

\begin{corollary}\label{non-tube->tube}
  If $G$, $\widetilde G$ and $H''$ are as in Theorem~\ref{rigidity-higgs},
  then the number of connected components of $\cM_{\max}(G)$ 
is bounded by the number of connected components of 
$\cM_{\max}(\widetilde G)$. In particular, if 
$\cM_{\max}(\widetilde G)$ is connected so is $\cM_{\max}(G)$.
\end{corollary}

\section{Morse theory on the  moduli space
of $G$-Higgs bundles}
\label{sec:morse-theory-moduli}

\subsection{Deformation theory}
\label{sec:deformation-theory}

Below we shall be doing analysis, in the form of Morse theory, on the
moduli spaces of $G$-Higgs bundles and therefore we need a description
of their tangent spaces.  This can be conveniently done using
hypercohomology of certain complexes of sheaves.  This idea probably
has its origin in Welters \cite{welters:1983}.  A convenient reference
is Biswas and Ramanan \cite{biswas-ramanan:1994}.

Let $(E,\varphi)$ be a $G$-Higgs bundle.  The \emph{deformation complex}
of $(E,\varphi)$ is the following complex of sheaves:
\begin{equation}\label{eq:def-complex}
  C^{\bullet}(E,\varphi)\colon E(\liehc) \xrightarrow{\ad(\varphi)}
  E(\liemc)\otimes K.
\end{equation}
Note that this makes sense because $[\lie{m}^{\C},\lie{h}^{\C}]
\subseteq \lie{m}^{\C}$.

The following result generalizes the fact that the infinitesimal
deformation space of a holomorphic vector bundle $V$ is
isomorphic to $H^1(\End V)$.

\begin{proposition}
  \label{prop:deform}
  The space of infinitesimal deformations of a $G$-Higgs bundle
  $(E,\varphi)$ is isomorphic to the hypercohomology group
  $\HH^1(C^{\bullet}(E,\varphi))$.
\end{proposition}

In particular, if $(E,\varphi)$ represents a non-singular point of the
moduli space $\mathcal{M}_d(G)$ then the tangent space at this point
is canonically isomorphic to $\HH^1(C^{\bullet}(E,\varphi))$.

For usual holomorphic vector bundles, the analogue of the following
result is the fact that the only endomorphisms of a stable bundle
are the constant multiples of the identity.

\begin{proposition}\label{prop:stable-implies-simple}
  Let $(E,\varphi)$ be a stable $G$-Higgs bundle which represents a
  smooth point of $\mathcal{M}_d(G)$.  Then
  \begin{displaymath}
    \HH^0(C^{\bullet}(E,\varphi)) = \HH^2(C^{\bullet}(E,\varphi))
    = 0.
  \end{displaymath}
\end{proposition}

\subsection{Morse theory}
\label{sec:morse-theory}

The idea of applying Morse theory to the study of moduli of
holomorphic vector bundles  has its origin in the fundamental
work of Atiyah and Bott \cite{atiyah-bott:1982}.  Here the moduli
space of stable bundles was studied using equivariant Morse theory on
the infinite dimensional space of unitary connections.  The use of
Morse theory in moduli spaces of Higgs bundles was introduced by
Hitchin \cite{hitchin:1987a}.  In this section we explain how to apply
these methods in moduli spaces of $G$-Higgs bundles.  In particular,
we give a criterion (Corollary~\ref{cor:proper}) for finding the local
minima of the Morse function, which is extremely useful in the context
of problem of counting connected components of the moduli space.

In order to define the Morse function, we shall consider the moduli
space $\mathcal{M}_d(G)$ of $G$-Higgs bundles from the gauge theory
point of view, as explained in Sec.~\ref{sec:correspondence}.  Thus we
identify $\mathcal{M}_d(G)$ with the moduli space of solutions to
Hitchin's equations (\ref{hitchin}).  From this point
of view it makes sense to define
\begin{equation}
  \label{eq:def-morse}
  \begin{aligned}
    f\colon\mathcal{M}_d(G) &\lra  \R,\\
    (d_A,\varphi) &\mapsto \norm{\varphi}^2,
  \end{aligned}
\end{equation}
where $\norm{\varphi}^2 = \int_X \abs{\varphi}^2 d\vol$ is the $L^2$-norm of
$\varphi$.  Note that this norm is well defined because $\abs{\varphi}^2$ is
invariant under $H$-gauge transformations.

The maps $f$ has its origin in symplectic geometry: away from the
singular locus of $\mathcal{M}_d(G)$ it is a moment map for the
hamiltonian $S^1$-action given by
\begin{displaymath}
  e^{i\theta}\colon(d_A,\varphi) \mapsto (d_A, e^{i\theta}\varphi).
\end{displaymath}
This fact is important for two reasons.  Firstly, a theorem of Frankel
\cite{frankel:1959} guarantees that, when $\mathcal{M}_d(G)$ is
smooth, $f$ is a perfect Bott--Morse function.  Secondly, the critical
points of $f$ are exactly the fixed points of the $S^1$-action.

The function $f$ can be used to obtain information about connected
components even when $\mathcal{M}_d(G)$ has singularities due to the
following result, proved by Hitchin \cite{hitchin:1987a}, using
Uhlenbeck's weak compactness theorem \cite{uhlenbeck:1982}.

\begin{proposition}\label{prop:proper}
  The function  $f\colon \mathcal{M}_d(G) \to \R$ is a proper map.
\end{proposition}

\begin{corollary}
  \label{cor:proper}
  Let $\mathcal{M} \subseteq \mathcal{M}_d(G)$ be a closed subspace
  and let $\mathcal{N} \subseteq \mathcal{M}$ be the subspace of local
  minima of $f$ on $\mathcal{M}$.  If $\mathcal{N}$ is connected, then
  so is $\mathcal{M}$.
\end{corollary}

\subsection{A criterion for minima}
\label{sec:criterion-minima}

In view of Corollary~\ref{cor:proper} it is clearly of fundamental
importance to have a criterion which allows one to identify the local
minima of the function $f$.  As has already been pointed out, the
critical points of $f$ are just the fixed points of the $S^1$-action
on $\mathcal{M}_d(G)$.  The $G$-Higgs bundles corresponding to fixed
points are the so-called \emph{Hodge bundles}
(see Hitchin \cite{hitchin:1987a,hitchin:1992} and Simpson
\cite{simpson:1992}), described in the following proposition.

\begin{proposition}
  \label{prop:vhs}
  A polystable $G$-Higgs bundle $(E,\varphi)$ corresponds to a fixed
  point of the action of $S^1$ on $\mathcal{M}_d(G)$ if and only if
  $(E,\varphi)$ is a \emph{Hodge bundle}, i.e.,
  there is a semi-simple element $\psi\in H^0(E(\lieh))$ and
  decompositions
  \begin{align*}
    E(\liehc) &= \bigoplus_k E(\liehc)_k,\\
    E(\liemc) &= \bigoplus_k E(\liemc)_k
  \end{align*}
  in eigen-bundles for $\psi$ such that
  \begin{displaymath}
    \psi_{|E(\liehc)_k} = ik\qquad\text{and}\qquad
    \psi_{|E(\liemc)_k} = ik,
  \end{displaymath}
  and, moreover, $[\psi,\varphi] = i\varphi$.
\end{proposition}

Notice that the condition $[\psi,\varphi] = i\varphi$ means that
\begin{displaymath}
  \varphi \in H^0(E(\liemc)_1 \otimes K).
\end{displaymath}
Hence, if $(E,\varphi)$ is a  Hodge bundle as
described in the preceding proposition, there is an induced
decomposition of the deformation complex $C^\bullet(E,\varphi$) defined
in \eqref{eq:def-complex}, as follows:
\begin{equation}\nonumber
  C^{\bullet}(E,\varphi) = \bigoplus_k C^{\bullet}_k(E,\varphi),
\end{equation}
where for each $k$ we define the complex
\begin{displaymath}
  C^{\bullet}_k(E,\varphi)\colon E(\liehc)_k
    \xrightarrow{\ad(\varphi)}E(\liemc)_{k+1}\otimes K.
\end{displaymath}
This decomposition gives us a corresponding decomposition of the
infinitesimal deformation space of $(E,\varphi)$:
\begin{equation}\nonumber
  \HH^1(C^{\bullet}(E,\varphi))
    = \bigoplus_k \HH^1(C^{\bullet}_k(E,\varphi)).
\end{equation}
This decomposition is important because of the following result.

\begin{proposition}
  Let $(E,\varphi)$ be a stable $G$-Higgs bundle which represents a
  non-singular point of $\mathcal{M}_d(G)$.  If $(E,\varphi)$ represents
  a fixed point of the $S^1$-action, then the eigenvalue $-k$
  eigenspace of the tangent space for the Hessian of $f$ is isomorphic
  to $\HH^1(C^{\bullet}_k(E,\varphi))$.  In particular, $(E,\varphi)$
  corresponds to a local minimum of $f$ if and only if
  \begin{displaymath}
    \HH^1(C^{\bullet}_k(E,\varphi)) = 0\ \forall k>0.
  \end{displaymath}
\end{proposition}

A key result proved in \cite{bradlow-garcia-prada-gothen:2003} gives a
very useful criterion for deciding when the hypercohomology
$\HH^1(C^{\bullet}_k(E,\varphi))$ vanishes.  In order to state this
result, it is convenient to use the Euler characteristic of the
complex $C^\bullet_k(E,\varphi)$ defined by
\begin{displaymath}
  \chi(C^\bullet_k(E,\varphi)) = \dim\HH^0(C^\bullet_k(E,\varphi)) -
  \dim\HH^1(C^\bullet_k(E,\varphi)) + \dim\HH^2(C^\bullet_k(E,\varphi)).
\end{displaymath}

\begin{theorem}[{\cite[Proposition~4.14]{bradlow-garcia-prada-gothen:2003}}]
  \label{thm:adjoint}
  Let $(E,\varphi)$ be a semi-stable $G$-Higgs bundle, fixed under the
  action of $S^1$.  Then
  \begin{displaymath}
    \chi(C^\bullet_k(E,\varphi)) \leq 0
  \end{displaymath}
  and equality holds if and only if $\ad(\varphi)\colon E(\liehc)_k \to
  E(\liemc)_{k+1}\otimes K$ is an isomorphism.
\end{theorem}

Together with Proposition~\ref{prop:stable-implies-simple},
Theorem~\ref{thm:adjoint} gives the promised criterion for finding the
local minima of $f$.

\begin{corollary}\label{cor:adjoint-minima}
  Let $(E,\varphi)$ be a stable $G$-Higgs bundle which represents a
  non-singular point of $\mathcal{M}_d(G)$.  Then $(E,\varphi)$
  represents a local minimum of $f$ if and only if
  $$
  \ad(\varphi)\colon E(\liehc)_k \lra E(\liemc)_{k+1}\otimes K
  $$
  is an isomorphism for all
  $k>0$.
\end{corollary}

\section{Minima and counting of components}
\label{sec:minima-counting}

\subsection{Minima}
\label{sec:minima}

It is the purpose of this section to describe the local minima of the
function $f\colon \mathcal{M}_d(G) \to \R$ defined in
\eqref{eq:def-morse} for these groups.

Recall from Remark~\ref{rem:positive-toledo} that there is no loss of
generality in assuming that the Toledo invariant is positive. Thus, to
all the results stated here for positive Toledo invariant, there are
parallel results for negative Toledo invariant.

Since it creates no extra difficulties, we state the following Theorem
for arbitrary (positive) Toledo invariant, even though we are
presently only interested in the case of maximal Toledo invariant.
The Theorem is proved by using the criterion given in
Corollary~\ref{cor:adjoint-minima}, together with an extra argument to
deal with strictly polystable $G$-Higgs bundles.

\begin{theorem}\label{thm:G-minima}
  Let $(E,\varphi=\beta+\gamma)$ be a polystable $G$-Higgs bundle with
  positive Toledo invariant.
  \begin{enumerate}
  \item\label{item:G-minima-1} If $G$ is one of the groups $\SU(p,q)$,
    $\SO^*(2n)$ and $\SO_0(2,n)$ (with $n \neq 2,3$) then
    $(E,\varphi=\beta+\gamma)$ represents a local minimum on
    $\mathcal{M}_d(G)$ if and only if $\beta=0$.
  \item\label{item:G-minima-2} Let $(E,\varphi=\beta+\gamma)$ be a polystable
    $\Sp(2n,\R)$-Higgs bundle with $n \neq 2$ and let
    $(\E,\Phi) = \left(V \oplus V^*,\left(
      \begin{smallmatrix}
        0 & \beta \\
        \gamma & 0
      \end{smallmatrix}\right)
    \right)$ be the associated Higgs vector bundle.  Then $(E,\varphi)$
    represents a local minimum of $f$ if and only if one of the
    following situations occurs:
    \begin{enumerate}
    \item\label{item:2a} The vanishing  $\beta=0$ holds.
    \item\label{item:2b} The number $n$ is odd and there is a square
      root $L$ of the canonical bundle $K$ and a decomposition $V =
      LK^{-2[n/2]} \oplus LK^{-2[n/2]+2} \oplus \dots \oplus
      LK^{2[n/2]}$ with respect to which
      \begin{equation}\label{eq:b}
        \gamma = \left(
        \begin{smallmatrix}
          0 & \cdots  & 1 \\
          \vdots & \iddots
          & \vdots \vspace{3pt}\\
          1 & \cdots & 0
        \end{smallmatrix} \right)
        \qquad\text{and}\qquad
        \beta = \left(
        \begin{smallmatrix}
          0 & \cdots  & 1 & 0\\
          \vdots & \iddots
          & & \vdots \vspace{-2pt} \\
          1 & &  & \vdots \vspace{3pt} \\
          0 & \cdots & \cdots  & 0
        \end{smallmatrix} \right).
    \end{equation}
    In this case, necessarily the Toledo invariant is maximal, i.e.\
    $\deg(V) = n(g-1)$.
  \item\label{item:2c} The number $n$ is even and there is a square
    root $L$ of the canonical bundle $K$ and a decomposition $V =
    L^{-1}K^{2-n} \oplus L^{-1}K^{4-n} \oplus \dots \oplus
    L^{-1}K^{n}$ with respect to which $\beta$ and $\gamma$ are given
    by \eqref{eq:b}.  Also in this case, necessarily we are in the
    situation of maximal Toledo invariant, $\deg(V) = n(g-1)$.
    \end{enumerate}
  \end{enumerate}
\end{theorem}

\begin{remark}
  The case of $G=\SU(1,1) \cong \Sp(2,\R) \cong \SL(2,\R)$ is covered
  by both (\ref{item:G-minima-1}) and (\ref{item:G-minima-2}) of
  Theorem~\ref{thm:G-minima}.  This case (together with $G=\SO_0(2,1)
  \cong \PSL(2,\R)$) was studied by Hitchin~\cite{hitchin:1987a}.
\end{remark}

\begin{remark}
  Recall from Remark~\ref{rem:so22} that we have excluded the group
  $\SO_0(2,2)$ from our considerations --- in fact, the results for
  this group do not fit into the general statement given in
  Theorem~\ref{thm:G-minima}.
\end{remark}

\begin{remark}\label{rem:teichmuller-sp2n}
  The minima for the split real group $\Sp(2n,\R)$ described in
  (\ref{item:2b}) and (\ref{item:2c}) of Theorem~\ref{thm:G-minima}
  are exactly the ones that belong to the Teichm\"uller components
  defined by Hitchin \cite{hitchin:1992}.
\end{remark}

It remains to deal with the special cases $\Sp(4,\R)$ and $\SO_0(2,3)$
(cf.\ Sec.~\ref{sec:sp-so-coincidence}).

\begin{theorem}\label{thm:sp4-so-minima}
  Let $G$ be one of the groups $\Sp(4,\R)$ or $\SO_0(2,3)$ and let
  $(E,\varphi=\beta+\gamma)$ be a polystable $G$-Higgs bundle with
  positive Toledo invariant.  Then $(E,\varphi)$ represents a local
  minimum of $f$ if and only if one of the following situations
  occurs:
    \begin{enumerate}
    \item The vanishing  $\beta=0$ holds.
    \item\label{item:so}
    If $G=\SO_0(2,3)$ and the associated Higgs vector bundle is
    $(V \oplus W,\Phi)$, then there are decompositions in line bundles
      \begin{displaymath}
        V = K \oplus K^{-1} \qquad\text{and}\qquad
        W = M \oplus \mathcal{O} \oplus M^{-1},
      \end{displaymath}
      with $0 < \deg(M) \leq 4g-4$.
      With respect to these decompositions,
      \begin{displaymath}
        Q_V = \left(
        \begin{smallmatrix}
          0  & 1 \\
          1  & 0
        \end{smallmatrix}\right)\qquad\text{and}\qquad
        Q_W = \left(
        \begin{smallmatrix}
          0 & 0 & 1 \\
          0 & 1 & 0 \\
          1 & 0 & 0
        \end{smallmatrix}\right),
      \end{displaymath}
      $\gamma$ is the canonical section $1$ of\,
      $\Hom(\mathcal{O},K^{-1})\otimes K\cong \mathcal{O}$
      while $\beta$ is a non-zero section of\, $\Hom(M,K)\otimes K$,
      i.e.\
      \begin{displaymath}
        \beta\in H^0(M^{-1}K^2).
      \end{displaymath}
      In this case, necessarily the Toledo invariant is maximal, i.e.,
      $d = 2g-2$.
    \item\label{item:sp4}
      If $G=\Sp(4,\R)$ and the associated Higgs vector bundle is
    $(V \oplus V^*,\Phi)$, then there is a decomposition in line bundles
    \begin{displaymath}
      V = N \oplus N^{-1}K,
    \end{displaymath}
     with $g-1 < \deg(N) \leq 3g-3$. With respect to this decomposition,
     $\gamma\in H^0(S^2V^*\otimes K)$ is given by the tautological
     section
     \begin{displaymath}
       \gamma=\left(
         \begin{smallmatrix}
           0 & 1 \\
           1 & 0
         \end{smallmatrix}\right)
     \end{displaymath}
     and $\beta\in H^0(S^2V\otimes K)$ is a non-vanishing section of
     the form
     \begin{displaymath}
       \beta=\left(
         \begin{smallmatrix}
           0 & 0 \\
           0 & \widetilde{\beta}
         \end{smallmatrix}\right)\qquad
       \text{with}\qquad \widetilde{\beta}\in H^0(N^{-2}K^3)
     \end{displaymath}
     Also in this case, necessarily the Toledo invariant is maximal,
     i.e.,  $d = 2g-2$.
    \end{enumerate}
\end{theorem}

\begin{remark}\label{rem:so-sp-lift}
  Following through the correspondence between $\Sp(4,\R)$-Higgs
  bundles and $\SO_0(2,3)$-Higgs bundles described in
  Sec.~\ref{sec:sp-so-coincidence}, one sees that a minimum for
  $\Sp(4,\R)$ of the type given in (\ref{item:sp4}) of
  Theorem~\ref{thm:sp4-so-minima} gives rise to a minimum of the type
  given in (\ref{item:so}) of the theorem with $M = N^2K^{-1}$.
  Thus, in particular, $\deg(M)$ is even and, in fact, one can see
  that the second Stiefel--Whitney class $w_2(W,Q_W)$ is exactly the
  modulo $2$ reduction of the degree of $M$, thus confirming (in the
  case of these minima) that this
  is the obstruction to lifting to $\Sp(4,\R)$.
\end{remark}

\begin{remark}
  In the case $G=\SO_0(2,3)$, the Cayley correspondence of
  Theorem~\ref{cayley-correspondence} becomes particularly simple to
  describe for the $\SO_0(2,3)$-Higgs bundles which are local minima (of
  maximal Toledo invariant, of course): the kernel of $\gamma$ is $F =
  M \oplus M^{-1}$ (cf.\ \eqref{extension}) and the restriction of
  $Q_W$ to $F$ is clearly non-degenerate.  Furthermore, $L_0=L^{-1}K
  =K^{-1}K$ is trivial, $\theta''=0$ and $\theta'\colon F \to K^2$ is
  given by $\beta\in H^0(M^{-1}K^2)=H^0(\Hom(M,K^2))$.
\end{remark}

\begin{remark}\label{rem:teichmuller-sp4-so}
  The minima for $\SO_0(2,3)$ which belong to a Teichm\"uller
  component are of the type described in (\ref{item:so}) of
  Theorem~\ref{thm:sp4-so-minima} with $\deg(M) = 4g-4$.  Note that,
  since $\beta\neq 0$, this forces $M=K^{-2}$.  Thus there is a unique
  such minimum and this lifts to a minimum for $\Sp(4,\R)$ because
  $\deg(M)$ is even (cf.\ Remark~\ref{rem:so-sp-lift}).

  The minima for $\Sp(2n,\R)$ which belong to a Teichm\"uller
  component are of the type described in (\ref{item:sp4}) of the
  Theorem with $\deg(N) = 3g-3$.  Note that, since $\beta\neq 0$, this
  means that $N=K^{3/2}$.  Hence we see that there are $2^{2g}$ such
  minima, corresponding to the choices of the square root of the
  canonical bundle $K$.  Clearly each of these minima are lifts of the
  unique minimum for $\SO_0(2,3)$.
\end{remark}

\subsection{The counting of components}
\label{sec:counting-components}

In this section we give the count of the number of components of
$\mathcal{M}_{\max}(G)$.  
Using Theorems \ref{rigidity-supq} and \ref{rigidityso*2n},
we can reduce the problem to the
case where $G/H$ is of tube type,
i.e., when $G$ is one of the groups $\SU(n,n)$, $\Sp(2n,\R)$,
$\SO(2,n)$ or $\SO^\ast(2n)$ with $n$ even.  
The number of connected
components of $\mathcal{M}_{\max}(G)$ is given in
Table~\ref{tab:components}.

The general strategy  for counting the components of
$\mathcal{M}_{\max}(G)$ when $G/H$ is of tube type  is as follows:
\begin{enumerate}
\item\label{item:step1} Use the Cayley correspondence of
  Theorem~\ref{cayley-correspondence} to obtain extra topological
  invariants, via the identification of $\mathcal{M}_{\max}(G)$ with
  the moduli space of $K^2$-twisted $G'$-Higgs pairs.  The relevant
  topological invariants can be read off Table~\ref{tab:tube} as those
  of bundles whose structure group is the maximal compact $H'\subseteq
  G'$.  This provides a subdivision
  \begin{displaymath}
    \mathcal{M}_{\max}(G) \cong \mathcal{M}_{K^2}(G')
    =\bigcup_{c} \mathcal{M}_{K^2,c}(G'),
  \end{displaymath}
  according to the values of these topological invariants.
\item Use the results of Sec.~\ref{sec:minima} to identify the local
  minima of $f$ on each of the subspaces $\mathcal{M}_{K^2,c}(G')$
  defined in (\ref{item:step1}).
\item For each subspace $\mathcal{M}_{K^2,c}(G')$, determine whether
  the space of local minima is connected and non-empty: if this is the
  case, then $\mathcal{M}_{K^2,c}(G')$ is a connected component of
  $\mathcal{M}_{\max}$.
\item When the subspace of local minima of $f$ on
  $\mathcal{M}_{K^2,c}(G')$ is not connected, find its connected
  components. It turns out that this only happens due to the presence
  of one or more Teichm\"uller components --- in this case,
  non-Teichm\"uller components with the same invariants may or may not
  exist.
\end{enumerate}

In the following we outline how this strategy is carried out for each
of the groups mentioned above.

\begin{itemize}

\item $G=\SU(n,n)$:

\end{itemize}

Consulting Table~\ref{tab:tube} we see that the relevant topological
invariants are those of $H'= \SU(n)\rtimes \Z_2$-bundles.  Since
$\pi_0(H')=\Z_2$ and $\pi_1(H')=\{1\}$, it follows that the invariant
takes values in $H^1(X;\Z_2) \cong \Z_2^{2g}$.  The corresponding
minima are given by (\ref{item:G-minima-1}) of
Theorem~\ref{thm:G-minima} as having $\beta=0$. In terms of the Cayley
correspondence \eqref{eq:cayley-su} this means that the pair
$(\widetilde{W},\widetilde{\theta})$ has $\widetilde{\theta}=0$.  Hence the
subspace of local minima is isomorphic to the moduli
space\footnote{Here, and in the following, one must check that the
  various stability conditions involved in defining the moduli spaces
  agree.  This is rather technical and we shall ignore this question
  in the present paper.} of vector bundles $\widetilde{W}$
such that $(\det \widetilde{W})^2=\mathcal{O}$.  This moduli
space has a connected component for each of the $2^{2g}$ choices of
square root of the trivial line bundle and it is not difficult to see
that this choice corresponds exactly to the value of the topological
invariant in $H^1(X;\Z_2) \cong \Z_2^{2g}$.  This gives a total of
$2^{2g}$ connected components of $\mathcal{M}_{\max}(G)$, as stated in
Table~\ref{tab:components}.  It is interesting to notice that this
analysis is also valid for the case $\SU(1,1) \cong \Sp(2,\R)$;
however it is somewhat special, in that each of the $2^{2g}$
components is a Teichm\"uller component (in fact isomorphic to
Teichm\"uller space).  This is obviously not the case when $n \neq 1$.

\begin{itemize}

\item  $G=\Sp(2n,\R)$:

\end{itemize}

In this case Table~\ref{tab:tube} shows that the relevant group is
$G'= \GL(n,\R)$.  Thus the invariants are the first and second
Stiefel--Whitney classes $w_1\in H^1(X;\Z_2) \cong \Z_2^{2g}$ and
$w_2\in H^2(X;\Z_2) \cong \Z_2$.  For each of the possible values of
$(w_1,w_2)$ there are minima of the type given in (\ref{item:2a}) of
Theorem~\ref{thm:G-minima}, i.e., with $\beta=0$.  Under the Cayley
correspondence \eqref{cayley} these correspond to pairs $(W,\theta)$
with $\theta=0$ (cf.\ Sec.~\ref{sec:cayley}) where $W$ is an
orthogonal bundle. Thus, for given $(w_1,w_2)$, the space of minima of
this type can be identified with the moduli space of
$\OO(n,\C)$-bundles with these invariants.  The moduli space of
principal bundles for a \emph{connected} group and fixed topological
type is known to be connected by Ramanathan
\cite[Proposition~4.2]{ramanathan:1975}.  Now, since $\OO(n,\C)$ is
not connected the result of Ramanathan cannot be applied directly,
however, all that is really required for his argument is that
semistability is an open condition and thus one obtains the desired
result (cf.\ \cite{gothen-oliveira:2005}).  Hence there are $2\cdot
2^{2g} = 2^{2g+1}$ connected components of
$\mathcal{M}_{\max}(\Sp(2n,\R))$ corresponding to minima with
$\beta=0$.

\emph{The case $n \geq 3$.} In this case there are additional minima
of the type described in (\ref{item:2b}) and (\ref{item:2c}) of
Theorem~\ref{thm:G-minima}.  Each of these minima are easily seen to
be the unique minimum of $f$ on a Teichm\"uller component (it can be
seen that they all have $w_2=0$, whereas $w_1=0$ when $n$ is even and
arbitrary for $n$ odd). Since there is such a minimum for each square
root of the canonical bundle, there are $2^{2g}$ Teichm\"uller
components.  Hence the total number of components is $3\cdot 2^{2g}$.

\emph{The case $n = 2$.} Here $G'=\GL(2,\R)$ and the maximal compact
subgroup is $H'=\OO(2)$.  One easily sees that the minima with
$\beta\neq 0$, described in (\ref{item:sp4}) of
Theorem~\ref{thm:sp4-so-minima} all have $w_1=0$.  Excluding the case
$w_1=0$ we thus have $(2^{2g}-1)\cdot 2 = 2^{2g+1}-2$ possible values
for the invariants, giving rise to the same number of components (the
subspaces of minima are connected, as above). When the first
Stiefel--Whitney class vanishes, there is a reduction of structure
group to $\SO(2) \cong S^1$ and then the second Stiefel--Whitney class
$w_2$ lifts to an integer invariant, namely the first Chern class
$c_1$ of the $S^1$-bundle.  For the case of the minima with $\beta=0$,
we have $c_1=0$ and again the space is connected. Hence, corresponding
to minima with $\beta=0$, there are in total $2^{2g+1}-1$ components.
For the minima with $\beta\neq 0$ described in (\ref{item:sp4}) of
Theorem~\ref{thm:sp4-so-minima} the value is $c_1=\deg(N)-(g-1)$.  It
follows from the bound $g-1 < \deg(N) \leq 3g-3$ that the
corresponding subspace is non-empty only for $c_1 \leq 2g-2$.  For
each $c_1$ satisfying $0 < c_1 < 2g-2$ it can be proved that the
subspace of minima is connected, thus showing that there is a unique
(non-Teichm\"uller, in fact) component for this value.  When $c_1 =
2g-2$ the minima belong to a Teichm\"uller component and there are
$2^{2g}$ of these (again depending on the choice of a square root of
the canonical bundle).  The total number of components is thus
$2^{2g+1}-1 + (2g-3) + 2^{2g} = 3\cdot 2^{2g} +2g - 4$.

\emph{The case $n = 1$.}  Since $\Sp(2,\R) \cong \SU(1,1)$ this case
has been treated above but, it can of course also be seen directly
that $\mathcal{M}_{\max}(\Sp(2,\R))$ has $2^{2g}$ Teichm\"uller-components.

\begin{itemize}

\item $G=\SO_0(2,n)$:

\end{itemize}

\emph{The case $n \geq 4$.}
In this case the group giving rise to the extra invariants is
$H'=\OO(n-1)$.  Thus the new invariants are again Stiefel--Whitney
classes $w_1$ and $w_2$.  From (\ref{item:G-minima-1}) of
Theorem~\ref{thm:G-minima} we know that the only minima are the ones
with $\beta=0$ so, exactly as explained for $\Sp(2n,\R)$ above, the
subspace of minima with given $(w_1,w_2)$ can be identified with the
moduli space of $\OO(n-1,\C)$-bundles, which is connected.  This gives
the total of $2\cdot 2^{2g} = 2^{2g+1}$ components of
$\mathcal{M}_{\max}(G)$.

\emph{The case $n = 3$.}
Here we have $G'= \SO(1,1)\times\SO(1,2)$ with maximal compact
subgroup $H'=\OO(2)$ and the analysis is parallel to the one given
above for $G=\Sp(4,\R)$.  The invariants are $(w_1,w_2)$ and all
minima with $w_1\neq = 0$ have $\beta=0$.  When $w_1=0$, the class
$w_2$ lifts to an integer class $c_1$ and for $c_1=0$ there is one
component whose minima have $\beta=0$.  The remaining components have
minima with $\beta\neq 0$ of the type described in (\ref{item:so}) of
Theorem~\ref{thm:sp4-so-minima} and the invariant is $c_1 = \deg(M)$.
For each allowed value $0<\deg(M)\leq 4g-4$ there is one connected component
(which lifts to $\Sp(4,\R)$ when $\deg(M)$ is even).   The
Teichm\"uller component occurs for $\deg(M)=4g-4$; in this case there
is just one, which is covered by the $2^{2g}$ projectively equivalent
components for $\Sp(4,\R)$.  Thus the total number of components is
$2^{2g+1}-1+4g-4 = 2^{2g+1}+4g-5$

\begin{remark}
The maximal compact subgroup of $\SO_0(2,n)$ is $H=\SO(2)\times
\SO(n)$.  Hence $\pi_1(H)\cong \Z \oplus \Z_2$ if $n\geq 3$. There
is thus a second invariant attached to an $\SO_0(2,n)$-Higgs bundle
(an element in $\Z_2$) in addition to the Toledo invariant. This can
be identified as the second Stiefel--Whitney class of the
$\SO(n,\C)$-bundle $(W, Q_W)$ in Table~\ref{tab:higgs-HSS}. We could
have used this invariant to distinguish two disjoint closed
subspaces in $\mathcal{M}_{\max}(\SO_0(2,n))$, but this is taken
into account by the Cayley correspondence.

\end{remark}

\emph{The case $n = 1$.}
In this case $G=\SO_0(2,1) \cong \PSL(2,\R) \cong \mathrm{PSp}(2,\R)$ and
there is a unique minimum with $\beta=0$.  This belongs to the
unique Teichm\"uller component and lifts to the $2^{2g}$ projectively
equivalent Teichm\"uller components for $G=\SL(2,\R)$.

\begin{itemize}

\item $G=\SO^*(2n)$, with $n=2m$:

\end{itemize}

In this case Theorem~\ref{thm:G-minima} shows that all minima
have $\beta=0$ and so from the Cayley correspondence (cf.\
Sec.~\ref{sec:cayley}) we deduce that the subspace of minima on
$\mathcal{M}_{\max}(G)$ can be identified with the moduli space of
principal $H'^{\C}$-bundles.  Since $H'$ is the simply connected group
$\Sp(n)$ (see Table~\ref{tab:tube}), it follows from Ramanathan
\cite{ramanathan:1975} that the space of minima is connected, thus
showing that $\mathcal{M}_{\max}(G)$ is also connected.

We now study the case when $G/H$ is not of tube type.

\begin{itemize}

\item $G=\SU(p,q)$ with $p\neq q$:

\end{itemize}

Without loss of generality we assume that $p<q$. 
From Theorem \ref{rigidity-supq} we see that $\mathcal{M}_{\max}(G)$
is isomorphic to the moduli space $\mathcal{M}'$ of $G$-Higgs bundles
of the form \eqref{eq:rigidsupqhiggs} with $\deg(V) = p (g-1)$ and
\begin{displaymath}
  \det(V)\otimes \det(W') \otimes(W'') \cong \mathcal{O}, 
\end{displaymath}
which means, as already noted, that every maximal $\SU(p,q)$-Higgs bundle
reduces to an $\SSS(\U(p,p) \times \U(q-p))$-Higgs bundle.

Define a map from $\mathcal{M}'$ to the Jacobian parameterizing line
bundles of degree $p(g-1)$ by
\begin{align*}
  \mathcal{M}' &\to J^{p(g-1)}(X)\\
  (V,W',W'',\beta,\gamma) &\mapsto \det(V)
\end{align*}
The fibre of this map over a line bundle $L$ is a product of the
moduli space 
$$M_{L^{-2}K^p}(\GL(q-p,\C))$$ 
of polystable bundles $W''$
of rank $q-p$ and fixed determinant $L^{-2}K^p$ and the moduli space
\begin{displaymath}
  \widetilde{\mathcal{M}} = \{(V,W',\beta,\gamma) \;:\; \det(V) = L \}
\end{displaymath}
which is a subspace of the moduli space
$\mathcal{M}_{L^2K^{-p}}(\U(p,p))$ of (maximal) $\U(p,p)$-Higgs
bundles with fixed determinant $L^2K^{-p}$.  Clearly
$\mathcal{M}_{L^2K^{-p}}(\U(p,p))$ is isomorphic to the moduli space
$\mathcal{M}_{\max}(\SU(p,p))$ and, similarly to this case,
$\mathcal{M}_{L^2K^{-p}}(\U(p,p))$ has $2^{2g}$ connected components,
indexed by the square roots of $L^2K^{-p}$.  Fixing $\det(V)$
obviously amounts to fixing one of these square roots and hence
$\widetilde{\mathcal{M}}$ is just one of these connected components.
The map to the Jacobian defined above is surjective since, if we are
given another line bundle $L\otimes L_0$, the map
\begin{displaymath}
  (V,W',W'',\beta,\gamma) \mapsto 
(V\otimes L_0^{1/p},W'\otimes L_0^{1/p},W''\otimes L_0^{-2/(q-p)},\beta,\gamma) 
\end{displaymath}
gives an isomorphism between the fibre over $L$ and the fibre over $L
\otimes L_0$.  Hence, the connectedness of the Jacobian, of
$M_{L^{-2}K^p}(\GL(q-p,\C))$ and  of $\widetilde{\mathcal{M}}$ imply that
$\mathcal{M}'$ is connected and, consequently,
$\mathcal{M}_{\max}(G)$ is connected.

\begin{itemize}

\item $G=\SO^*(2n)$, with $n=2m+1$:

\end{itemize}

From Theorem \ref{rigidityso*2n} and the fact that the Jacobian
$J(X)$ is connected, we have that $\mathcal{M}(\SO^*(4m+2))$ is 
connected since, as we have seen,  $\mathcal{M}(\SO^*(4m))$ 
is connected.

\newpage
\section{ Tables}
\label{sec:tables}

\begin{table}[htbp]
\begin{tabular}{|c|c|c|c|}
\hline\raisebox{-8pt}{}
$G$ & $H$ &  $H^\C$ & $\liem^\C=\liemp+ \liemm$ \\
\hline\hline\raisebox{-8pt}{}
$\SU(p,q)$ & $\SSS(\U(p)\times \U(q))$ & $\SSS(\GL(p,\C)\times \GL(q,\C))$ &
$\Hom(\C^q,\C^p)+ \Hom(\C^p,\C^q)$ \\
\hline\raisebox{-8pt}{}
$\Sp(2n,\R)$ &  $\U(n)$ & $\GL(n,\C)$ &
$S^2(\C^n) + S^2({\C^n}^*)$ \\
\hline\raisebox{-8pt}{}
$\SO^*(2n)$ &  $\U(n)$  & $\GL(n,\C)$ &
$\Lambda^2(\C^n) + \Lambda^2({\C^n}^*)$ \\
\hline\raisebox{-8pt}{}
$\SO_0(2,n)$ & $\SO(2)\times \SO(n)$ & $\SO(2,\C)\times \SO(n,\C)$ &
 $\Hom(\C^n,\C)+ \Hom(\C,\C^n)$ \\
\hline

\end{tabular}
\vspace{12pt}

\caption{Irreducible classical Hermitian symmetric spaces $G/H$}
    \label{tab:HSS}
\end{table}

\begin{table}[htbp]
\begin{tabular}{|>$c<$||>$c<$|c|c|}
  \hline\raisebox{-20pt}{}
  G & \#\pi_0(\mathcal{M}_{\mathrm{max}}(G)) &
    \parbox[t]{61pt}{\centering Teichm\"uller components}&  Reference\\
  \hline\hline\raisebox{-20pt}{}
  \SU(n,n) & 2^{2g}
     & \parbox[t]{94pt}{\centering -- \\ {\footnotesize ($2^{2g}$ if $n=1$)}}
     & \cite{bradlow-garcia-prada-gothen:2003,markman-xia:2002} \\
  \hline\raisebox{-8pt}{}
  \SU(p,q) \quad (p \neq q)  & 1
     & --
     & \cite{bradlow-garcia-prada-gothen:2003,xia:2000,xia:2003} \\
  \hline\raisebox{-8pt}{}
  \Sp(2n,\R) \quad (n\geq 3)
  & 3\cdot 2^{2g}
    & $2^{2g}$
    & \cite{garcia-prada-gothen-mundet:2005} \\
  \hline\raisebox{-8pt}{}
   \SO_0(2,n) \quad (n\geq 4)
    & 2^{2g+1} & -- & \cite{bradlow-garcia-prada-gothen:preparation-1} \\
  \hline\raisebox{-8pt}{}
  \SO^*(2n) & 1 & -- & \cite{bradlow-garcia-prada-gothen:preparation-1} \\
  \hline\raisebox{-8pt}{}
  \Sp(4,\R) \cong \Spin_0(2,3)
    & 3\cdot 2^{2g} + 2g-4
    & $2^{2g}$
    & \cite{gothen:2001} \\
  \hline\raisebox{-8pt}{}
  \SO_0(2,3)
    & 2^{2g+1} + 4g-5
    & $1$
    & \cite{bradlow-garcia-prada-gothen:preparation-1} \\
  \hline\raisebox{-8pt}{}
  \Sp(2,\R) \cong \SL(2,\R)
    & 2^{2g}
    & $2^{2g}$
    & \cite{goldman:1988,hitchin:1987a} \\
  \hline\raisebox{-8pt}{}
  \SO_0(2,1) \cong \PSL(2,\R)
    & 1
    & $1$
    & \cite{goldman:1988,hitchin:1987a} \\
  \hline
\end{tabular}
\vspace{12pt}

\caption{Components of $\mathcal{M}_{\max}(G)$}
    \label{tab:components}

\end{table}

\begin{landscape}
\begin{table}[htbp]

\begin{tabular}{|c||c|c|c|c|}

\hline
\raisebox{-8pt}{}


$G$

&

$\SU(p,q)$

&

$\Sp(2n,\R)$

&

$\SO^*(2n)$

&

$\SO_0(2,n)$

\\

\hline\hline

$(E,\varphi)$

&

$V$: rank $p$ bundle

&

$V$: rank $n$ bundle

&

$V$: rank $n$ bundle

&

$
\big( V=L\oplus L^{-1},Q_V=
\left(\begin{smallmatrix}{} 0 & 1\\ 1 & 0 \end{smallmatrix}\right) \big)
$

\\

& $W$: rank $q$ bundle

&

&

&

\parbox[t]{93pt}{\centering $(W,Q_W)$: rank $n$ orthogonal bundle}

\\

&

$\det  V\otimes \det  W=\cO$

&

&

&

$L$: line bundle;  $\det  W=\cO$

\\

$\varphi=\beta+\gamma$

&

$\beta\in H^0(\Hom(W,V)\otimes K)$

&

$\beta\in H^0(S^2 V\otimes K)$

&

$\beta\in H^0(\Lambda^2 V\otimes K)$

&

$\beta\in H^0(\Hom(W,L)\otimes K)$

\\

&

$\gamma\in H^0(\Hom(V,W)\otimes K)$

&

$\gamma\in H^0(S^2 V^*\otimes K)$

&

$\gamma\in H^0(\Lambda^2 V^*\otimes K)$

&

$\gamma\in H^0(\Hom(W,L^{-1})\otimes K)$

\raisebox{-8pt}{}

\\

\hline


$G^\C\subset \SL(N,\C)$

&

&

&

&

\\

$\E=E(\C^N)$

&

$\E= V\oplus W$

&

$\E= V\oplus V^*$

&

$\E= V\oplus V^*$

&

$\E= V\oplus W$

\\

$\Phi\in H^0(\End \E \otimes K)$

&

$\Phi = \left(\begin{smallmatrix}{} 0 & \beta \\ \gamma  & 0 \end{smallmatrix}\right)$

&

$\Phi = \left(\begin{smallmatrix}{} 0 & \beta \\ \gamma  & 0 \end{smallmatrix}\right)$

&

$\Phi = \left(\begin{smallmatrix}{} 0 & \beta \\ \gamma  & 0 \end{smallmatrix}\right)$

&

$\Phi = \left(\begin{smallmatrix}{} 0 & 0 & \beta \\ 0 & 0 &  \gamma
\\ -\gamma^t & -\beta^t & 0
 \end{smallmatrix}\right)$

\raisebox{-18pt}{}

\\

\hline


Toledo  invariant

&

$d=\deg V=-\deg W$

&

$d=\deg V$

&

$d=\deg V$

&

$d=\deg L$

\raisebox{-8pt}{}

\\

\hline

\parbox[t]{71pt}{\centering Milnor--Wood inequality}

&

&

&

&

\\

$|d|\leq d_{\max}$

&

$|d|\leq \min\{p,q\}(g-1)$

&

$|d|\leq n(g-1)$

&

$|d|\leq [\frac{n}{2}](2g-2)$

&

$|d|\leq 2g-2$

\raisebox{-8pt}{}
\\

\hline

\end{tabular}

\vspace{12pt}

\caption{Higgs bundles for irreducible classical symmetric spaces $G/H$}
    \label{tab:higgs-HSS}
\end{table}

\begin{table}[htbp]
\begin{tabular}{|c|c|c|c|c|c|c|}
\hline\raisebox{-8pt}{}
$G$ & $H$ &  $G'$ & $H'$ & $\Sh=H/H'$ & $\liem'$ &$ \liem'^\C$\\
\hline\hline\raisebox{-8pt}{}
$\SU(n,n)$ & $\SSS(\U(n)\times \U(n))$ & 
\parbox[t]{110pt}{$\{A \in \GL(n,\C) \suchthat$ \\ \mbox{}\hfill $\det(A)^2 = 1 \}$ \\
  $ \cong \SL(n,\C) \rtimes \Z_2$} & 
\parbox[t]{80pt}{$\{A \in \U(n) \suchthat$ \\ \mbox{}\hfill $\det(A)^2 = 1 \}$ \\
$ \cong \SU(n) \rtimes \Z_2$} &
$\U(n)$  & $\Herm(n,\C)$  & $\Mat(n,\C)$  \\
\hline\raisebox{-8pt}{}
$\Sp(2n,\R)$ &  $\U(n)$  & $\GL(n,\R)$ & $\OO(n)$ &
$\U(n)/\OO(n)$  & $\Sym(n,\R)$ &  $\Sym(n,\C)$  \\
\hline\raisebox{-8pt}{}
$\SO^*(2n)\;, n=2m$ &  $\U(n)$  & $\U^*(n)$ & $\Sp(n)$ &
$\U(n)/\Sp(n)$  & $\Herm(m,\HH)$ &  $\Skew(n,\C)$  \\
\hline\raisebox{-14pt}{}
$\SO_0(2,n)$ & $\SO(2)\times \SO(n)$ & $\SO_0(1,1)\times \SO(1,n-1)$ & $\OO(n-1)$ &
\parbox[b]{80pt}{\mbox{}\\$\dfrac{\U(1)\times S^{n-1}}{\Z_2}$}  & $\R\times \R^{n-1}$ &  $\C\times \C^{n-1}$  \\
\hline

\end{tabular}

\vspace{12pt}

\caption{Irreducible classical tube type Hermitian symmetric spaces $G/H$}
    \label{tab:tube}
\end{table}

\begin{table}[htbp]
\begin{tabular}{|c|c|c|c|c|c|c|}
\hline\raisebox{-8pt}{}
$G$ & $H$ &  $H'$ & $\widetilde G$ & $\widetilde H$ & $\widetilde H'$ &$ H''= H'/\widetilde H'$\\
\hline\hline\raisebox{-8pt}{}
$\SU(p,q)\;, p<q $ & $\SSS(\U(p)\times \U(q))$ & 
\parbox[t]{135pt}{$\{(A,B) \in \U(p)\times \U(q-p):$ \\ \mbox{}\hfill $\det(A)^2\det(B) = 1 \}$ \\
  $ \cong \SSS(\U(p)\times\U(q-p))\rtimes \Z_2$} 
& $\SU(p,p)$ &  $\SSS(\U(p)\times \U(p))$  & $\SU(p)\rtimes \Z_2$  
& $\SSS(\U(1)\times\U(q-p))$  \\
\hline\raisebox{-8pt}{}
$\SO^*(4m+2)$ &  $\U(2m+1)$  & $\Sp(2m)\times \U(1)$ &
$\SO^*(4m) $ &  $\U(2m)$  & $\Sp(2m)$ &  $\U(1)$  \\
\hline
\end{tabular}
\vspace{12pt}

\caption{Irreducible classical non-tube type Hermitian symmetric spaces $G/H$}
    \label{tab:non-tube}
\end{table}

\end{landscape}


\providecommand{\bysame}{\leavevmode\hbox to3em{\hrulefill}\thinspace}

\end{document}